\documentclass[11pt]{amsart}

\pdfoutput=1

\usepackage[T1]{fontenc}
\usepackage{lmodern}
\usepackage[utf8x]{inputenc}
\usepackage{graphicx}
\DeclareGraphicsExtensions{.pdf,.png,.jpg}
\usepackage{latexsym, stmaryrd}
\usepackage{verbatim}
\usepackage{amsmath}
\usepackage{amsthm}
\usepackage{amssymb}
\usepackage{hyperref}
\usepackage{pstricks-add}
\usepackage{enumitem}
\usepackage{bm}

\numberwithin{equation}{section}
\newtheorem{theorem}{Theorem}[section]

\newtheorem{remark}[theorem]{Remark}

\DeclareMathOperator{\diverg}{div}

\newcommand{\ve}{\varepsilon}

\newcommand{\ie}{\emph{i.e. }}
\newcommand{\Q}{\mathcal{Q}}
\newcommand{\G}{\mathcal{G}}
\newcommand{\Tau}{\mathcal{T}}
\newcommand{\RR}{\mathbb{R}}

\newcommand{\TT}{\mathbb{T}}
\newcommand{\SSS}{\mathbb{S}}

\usepackage{geometry}
\author[F. Filbet]{Francis Filbet}
\address{Francis Filbet \\
Universit\'e de Lyon \\
CNRS UMR 5208 \\
Universit\'e Lyon 1 \\
Institut Camille Jordan \\
43 blvd. du 11 novembre 1918 \\
F-69622 Villeurbanne cedex \\
France.}
\email{filbet@math.univ-lyon1.fr}

\author[Th. Rey]{Thomas Rey}
\address{Thomas Rey \\
Universit\'e de Lyon \\
CNRS UMR 5208 \\
Universit\'e Lyon 1 \\
Institut Camille Jordan \\
43 blvd. du 11 novembre 1918 \\
F-69622 Villeurbanne cedex \\
France.}
\email{rey@math.univ-lyon1.fr}

\title[A Rescaling Velocity Method for Dissipative Kinetic Equations]{A Rescaling Velocity Method for Dissipative Kinetic Equations \\ Applications to Granular Media}

	

\date{}%{Preliminary version of \today}

\begin{document}

%%%%%%%%%%%%%%%%%%%%%%%%%%%%%%%%%%%%%%%%%%%%%%%%%%%%%%%%%%%%%%%%%%%%%%%%%%%%%%%%%%%%%%%%%%%%%%%%%%%%%%%%%%%%%%%%%%%%	
\begin{abstract}
We present a new numerical algorithm based on a relative energy scaling for collisional kinetic equations allowing to study numerically their long time behavior, without the usual problems related to the change of scales in velocity variables. 
It is based on the knowledge of the hydrodynamic limit of the model considered, but is able to compute solutions for either dilute or dense regimes.
Several applications are presented for Boltzmann-like equations. This method is particularly efficient for  numerical simulations of the granular gases equation with dissipative energy: it allows to study accurately the long time behavior of this equation and is very well suited for the study of clustering phenomena.
\end{abstract}	
%%%%%%%%%%%%%%%%%%%%%%%%%%%%%%%%%%%%%%%%%%%%%%%%%%%%%%%%%%%%%%%%%%%%%%%%%%%%%%%%%%%%%%%%%%%%%%%%%%%%%%%%%%%%%%%%%%%%	
\maketitle
%%%%%%%%%%%%%%%%%%%%%%%%%%%%%%%%%%%%%%%%%%%%%%%%%%%%%%%%%%%%%%%%%%%%%%%%%%%%%%%%%%%%%%%%%%%%%%%%%%%%%%%%%%%%%%%%%%%%	
\tableofcontents
%%%%%%%%%%%%%%%%%%%%%%%%%%%%%%%%%%%%%%%%%%%%%%%%%%%%%%%%%%%%%%%%%%%%%%%%%%%%%%%%%%%%%%%%%%%%%%%%%%%%%%%%%%%%%%%%%%%%

\section{Introduction}
We are interested in numerical simulations of the long time behavior of collisional kinetic equations such as the Boltzmann equation for granular gases. To this end, we introduce a new technique which is based on the information provided by the hydrodynamic fields computed from a macroscopic model corresponding to the original kinetic equation. Then we simply rescale the kinetic equation according to these hydrodynamic quantities. The reason to do so is that the change of scales in velocity is a challenging numerical problem when one wants to deal with solutions to dissipative kinetic equations on a fixed grid. Indeed, most of the usual deterministic methods fails to capture the correct long time behavior because of concentration or spreading over the velocity space.

Recently, F. Filbet and G. Russo proposed in~\cite{Filbet:2004} a rescaling method for space homogeneous kinetic equations on a fixed grid. This idea is mainly based on the self-similar behavior of the solution to the kinetic equation. However for space non homogeneous case, the situation is much more complicated and this method cannot be applied since the transport operator and the boundary conditions break down this self-similar behavior. Here we propose a very simple idea based on a rescaling of the kinetic equation according to its hydrodynamic limit.   We will see that this approach can actually be applied to various type of collisional kinetic equations such as elastic and inelastic Boltzmann equation (also known as the \emph{granular gases equation}).

Among popular methods for numerical simulations of collisional gases, the most widely used are meshless. On the one hand, the Molecular Dynamics algorithm is a deterministic method up to the random initialization of the particles which is valid in a dense regime \cite{haile:1992}. This method has been used to describe the apparition of shocks in supersonic sand and reproduces very accurately experimentation results \cite{Rericha:2002}. On the other hand, Direct Simulation Monte Carlo (DSMC) methods are stochastic algorithms working either in dense or rarefied regimes. We refer for instance to \cite{Bird:1994} for a complete review of this topic. 
These methods are really efficient in term of computational cost since their complexity grows linearly with the number $N$ of particles but they are rather inaccurate since the order of accuracy is about $\mathcal{O}(1/\sqrt{N})$.  Another approach consists in a direct resolution of the Boltzmann operator on a phase space grid. For instance, deterministic and highly accurate methods based on a spectral discretization of the collisional operator have been proposed by F. Filbet, G. Naldi, L. Pareschi, G. Toscani and G. Russo in \cite{PaRu:SINUM:2000,Naldi:2003, filbet:2005} for the space-homogeneous setting. Although being of complexity $ \mathcal{O}(N^2)$, they are spectrally accurate and then need very few points to be precise. Another spectral method, inspired of the direct fast Fourier transform approach of A.V. Bobylev and S. Rjasanow~\cite{Bobylev:1997,BoRj:98,BoRj:99}  using Lagrange multiplier to improve the conservations was also used for space-homogeneous simulations in 3D by I. Gamba and S.H. Tharkabhushanam in~\cite{Gamba:2009}.

Due to the large number of particles involved in the study of rarefied gas dynamics, we shall adopt a statistical physics' point of view, by the use of Boltzmann-like kinetic equations.  For a given nonnegative initial condition $f_0$, we will consider a particle distribution function $f^\ve = f^\ve(t,x,v)$, for $t \geq 0$, $x \in \Omega \subset \mathbb{R}^{d_x}$ and $v \in \mathbb{R}^{d_v}$, solution to the initial-boundary value problem
\begin{equation} 
\label{eqCollision}
\left\{ \begin{array}{l}
\displaystyle\frac{\partial f^\ve}{\partial t} + v \cdot \nabla_x f^\ve \,=\, \frac{1}{\ve}\;\Q(f^\ve, f^\ve), 
\\
\;
\\
 f^\ve(0, x, v) = f_{0}(x,v).
\end{array}\right.
\end{equation}
where the collision operator $\Q$ is a  Boltzmann-like operator, which preserves at least mass and momentum. This equation  describes numerous models such as the Boltzmann equation for elastic and inelastic collisions or Fokker-Planck-Landau type equations.  The parameter $\ve > 0$ is the dimensionless Knudsen number, that is the ratio between the mean free path of particles before a collision and the length scale of observation. It measures the rarefaction of the gas: the gas is \emph{rarefied} if $\ve \sim 1$ and \emph{dense} if $\ve \ll 1$.  The open  set $\Omega$ is a bounded Lipschitz-continuous domain of $\RR^{d_x}$, which means  that the model (\ref{eqCollision}) has to be supplemented with boundary conditions described later.

The article is organized as follows. We present in Section \ref{secRescMeth} the rescaling velocity method and give two possible choices of scaling functions, the coupled Filbet-Russo scaling and the new uncoupled macroscopic one. We choose to use the latter in the rest of the article, and apply it to three different collision operators in Section \ref{secApp}: the Boltzmann operator, the full granular gases operator and a simplified BGK-like granular gases operator. 
Subsequently, we introduce in Section \ref{secBoundCond} the spatial boundary conditions we have to impose on $\Omega$: the so-called \emph{Maxwell} boundary conditions, a convex combination of specular and diffusive reflections at the boundary. We present these conditions for the kinetic model (in classical and rescaled variables) and for the macroscopic equations.
In Section \ref{secAlgo}, we introduce the numerical methods we shall use for the discretization of each term of the problem: the collision operator, the conservative transport term and the system of conservation laws with source term.
Finally, we present in Section \ref{secNumRes} some numerical results concerning the full and simplified granular gases models, in space homogeneous and inhomogeneous settings,  with different geometries.

\section{The Rescaling Velocity Method}
\label{secRescMeth}

This section is devoted to the presentation of a new  scaling for equation~\eqref{eqCollision} allowing to follow the change of scales in velocity. It is an extension to the space-dependent setting of the method first introduced in \cite{Filbet:2004}, using the relative kinetic energy as a scaling function. 
This method was reminiscent from the work of A. Bobylev, J.A. Carrillo, and I. Gamba  \cite{bobylev:2000} about Enskog-like inelastic interactions models. Indeed, in one section of this work, the authors scaled the solution of the space homogeneous collision equation by its thermal velocity, in order to study a drift-collision equation, where no blow-up occurs.
The same technique was also used by  S. Mischler and C. Mouhot in \cite{Mischler:20062} to prove the existence of self-similar solutions to the granular gases equation.
This temperature-dependent scaling has also been introduced by J.A. Carrillo and G. Toscani in~\cite{Carrillo:2000} and J.A. Carrillo, M. Di Francesco and G. Toscani in~\cite{Carrillo:2006Self} for nonlinear diffusion equations. 
It was finally extended by J.A. Carrillo and J. V\'azquez in \cite{Carrillo:2007a} to show the apparition of chaotic behavior for a precise (constructive) nonlinearity : the self-similar profile of this equation ``oscillates'' between Gaussian (heat equation) and Zel'dovich-Kompaneets-Barenblatt (porous medium equation) profiles.

Let us drop for simplicity the $\ve$-dependence in $f$. For a given positive function $\omega:\RR^+\times\Omega\mapsto \RR^+$, we introduce a new distribution $ g (t, x, \xi)$ by setting
\begin{equation} 
\label{defRescaling}
f (t,x,v) \,=\, \frac{1}{\omega^{d_v}} \,g \left (t,x,\frac{v}{\omega}\right ),
\end{equation}		
where the function $\omega$ is assumed to be an accurate measure of the ``support'' or  scale  of the distribution $f$ in velocity variables. Then according to this scaling, the distribution $g $ should naturally ``follow'' either the concentration or the spreading in velocity of the distribution $f$.

Let us now derive the kinetic equation verified by the distribution $g$.  Differentiating relation \eqref{defRescaling}  with respect to time yields
\begin{equation*}
\frac{\partial f}{\partial t} = \frac{1}{\omega^{d_v}}\left [ \frac{\partial g}{\partial t} - \frac{1}{\omega} \frac{\partial \omega}{\partial t} \diverg_\xi (\xi \, g)\right ]. 
\end{equation*}
Provided that $v = \omega \xi $, one also has
\begin{equation*}
v \cdot \nabla_x f = \frac{\xi }{\omega^{{d_v}-1}}\cdot\left [ \nabla_x g - \frac{1}{\omega} \nabla_x \omega \diverg_\xi (\xi \, g)\right ]. 
\end{equation*}
Then, if $f$ is solution to~\eqref{eqCollision}, the distribution $g$ given by~\eqref{defRescaling} is solution to the following equation:
\begin{equation*}
\frac{\partial g}{\partial t} + \omega \,\xi \cdot \nabla_x g  - \frac{1}{\omega} \left [ \left ( \frac{\partial \omega}{\partial t} + \omega \, \xi \cdot \nabla_x \omega \right ) \diverg_\xi(\xi \, g) \right ] = \frac{\omega^{d_v}}{\ve} \, \widetilde{\Q}(g,g),
\end{equation*}
where $ \widetilde \Q$ is such that $ \widetilde \Q(g,g) = \Q(f,f)$. This equation can actually be written in the more convenient conservative form:
\begin{equation} 
\label{eqGcons}
\frac{\partial g}{\partial t} + \diverg_x \left ( \omega \,\xi \, g\right )  - \diverg_\xi \left [ \left ( \frac{1}{\omega} \frac{\partial \omega}{\partial t} \, \xi + \xi \otimes \xi \, \nabla_x \omega \right ) g \right ] = \frac{\omega^{d_v}}{\ve} \, \widetilde{\Q}(g,g).
\end{equation}

In order to make this rescaling efficient, the main difficulty is now to choose an  appropriate scaling function $\omega$ to define completely the distribution $g$. 
			
\subsection{The Generalized Filbet-Russo Scaling}

The first and natural idea follows from the work~\cite{Filbet:2004}. It consists in computing the function $\omega$ directly from the distribution function $f$, by setting
\begin{equation} 
\label{eqWcoupled}
\omega \,\,:=\,\, \sqrt{\frac{ 2\,E_f}{{d_v}\,\rho_f}},
\end{equation}
where  $\rho_f$ and $E_f$ are respectively the local mass and kinetic energy of $f$ 
\begin{equation}
\label{def:rhoE}
 \rho_f \,:=\, \int_{\mathbb{R}^{d_v}} f(v) \, dv, \quad \ E_f \,:=\, \int_{\mathbb{R}^{d_v}} f(v) \, \frac{|v|^2}{2} \, dv
\end{equation}
We also define for the sequel the mean velocity field $\bm{u}_f$ and temperature $T_f$ by 
\begin{equation}
\label{def:uT}
 \bm{u}_f \,:=\, \frac{1}{\rho_f}\int_{\mathbb{R}^{d_v}} f(v) \, v\,dv, \quad T_f := \frac{1}{{d_v}\,\rho_f}\left( 2\,E_f \,-\, \rho_f \,|\bm{u}_f|^2 \right).  
\end{equation}
Assuming that $f$ is nonnegative, the quantity $\omega$ will then provide  correct information on its support: if $f$ is concentrated, $\omega$ will be small, whereas it will be large for scattered distributions. This approach has been shown to be very accurate for the space-homogeneous setting in~\cite{Filbet:2004}, but it is difficult to extend to our case because the definition of $\omega$ yields very restrictive constraints on the moments of $g$, namely
%
%Indeed using the change of variable $v \to \omega \, \xi$ in \eqref{def:rhoE}, we find that
%\begin{equation*}
%\omega^2(t,x) \, =\, \frac{2\,E_f(t,x)}{d\,\rho_f(t,x)} \,=\, \frac{2}{d}\,\frac{E_g(t,x)}{\rho_g(t,x)}\, \omega^2(t,x)
%\end{equation*}
%and then for a distribution function $g$  given by~\eqref{defRescaling}, the condition (\ref{eqWcoupled}) imposes that 
%\begin{equation*} 
%\rho_g(t,x) \,=\, \frac{2}{d}\, E_g(t,x),\quad \ \forall (t,x) \in \RR_+ \times \Omega.
%\end{equation*}
%Stated in another form, it means that
\begin{equation*}
  \int_{\RR^{d_v}} g(t,x,\xi) \, d\xi \,=\, \frac{1}{{d_v}} \int_{\RR^{d_v}} g(t,x,\xi) \, |\xi|^2 \, d\xi,\quad \ \forall (t,x) \in \RR_+ \times \Omega.
\end{equation*}
%This condition gives a very strong constraint on $g$ since it has to be satisfied locally both in time and space. From a numerical point of view, it means that the time and space discretizations have to be done in such a way that this condition is also  preserved at the discrete level, which becomes very restrictive on the choice of the numerical method. Moreover, given that $\omega$ depends on the distribution $g$ itself, the linear free transport equation satisfied by $f$ becomes a strongly nonlinear transport equation in $g$, which can be difficult to solve numerically.

\subsection{The Macroscopic Scaling}	
In order to relax this constraint on the distribution function $g$, we propose a different choice of the function $\omega$. Moreover, to avoid a strong coupling between $\omega$ and $g$ leading to nonlinear terms, we do not want to compute $\omega$ directly from the distribution function $g$ itself but only from macroscopic quantities which are assumed to be close enough to the one computed from $f$. 
Therefore, a good candidate will be a solution to the system of macroscopic equations obtained from a closure of the kinetic model \eqref{eqCollision}. But, as we want to develop a rescaling velocity method, what we really need is an accurate representation of the support of the initial condition, and not the good hydrodynamic limit of the model.

To find this system of equation, we choose to assume that, in the asymptotic limit $\ve \to 0$, the distribution function $f$ solution to (\ref{eqCollision}) converges towards a steady state $\mathcal{M}_{\rho, \bm{u}, T}$, where $\rho$ is the \emph{mass},  $ \bm{u}$ the \emph{velocity} and $T$ the \emph{temperature}. Then, a first order approximation with respect to $\ve$ would be
\begin{equation*}
f^{\ve} \,=\, \mathcal{M}_{\rho, \bm{u}, T} \,+\, \mathcal O (\ve),
\end{equation*}
where the macroscopic quantities are solutions to \cite{Golse:2005}
\begin{equation}
	\label{eqSysConsUncoupled}
	\left\{ \begin{aligned}
		& \partial_t \rho  + \diverg_x (\rho \, \bm{u} ) \,=\, 0, 
		\\
		& \,
		\\
		& \partial_t(\rho \, \bm{u} ) + \diverg_x  \left(\rho \, \bm{u} \otimes \bm{u}\,+\, p\,{\rm\bf I} \right)  \,=\, 0, 
		\\
		& \,
		\\
		& \partial_t E + \diverg_x \left(\bm{u}\, (  E \,+\, p ) \right) \,=\, G(\rho, \bm{u}, E).
	\end{aligned} \right.
\end{equation}
In this  last system, $E$ is the kinetic energy given by 
$$
2 E \,=\, {d_v} \rho \, T + \rho \, |\bm{u}|^2,
$$ 
 $p > 0$ is  the pressure and $G \in \RR$ the ``energy'' of the collision operator.
A natural choice for $\omega$ is to set
\begin{equation} 
\label{eqWclose}
\omega := \sqrt{\frac{2\,E}{{d_v}\,\rho}}.
\end{equation}
The scaling function $\omega$ becomes independent of the original distribution but still follows the change of scales in velocity, provided that the closure has been done properly.
Once more, we want to insist on the fact that this quantity does not have to be the exact ``temperature'' of the limit equation, but only an information on the support of the distribution $f$.  The conditions to apply the method are then to have some good estimates (by above and below if possible) of the change of scales of the equation, through the temperature of its solutions.
The system \eqref{eqGcons}-\eqref{eqSysConsUncoupled} is now uncoupled which is better both for a numerical and mathematical point of view. Moreover, it will be easier to solve thanks to the hyperbolic structure of~\eqref{eqSysConsUncoupled}.

\begin{remark}
  In the case where no precise estimates of the dissipation term $G$ are known, one possibility would be to close the moment equations with a Maxwellian distribution $\mathcal M_{\rho, \bm{u}, T}$, even if we can't compute explicitly the kinetic energy of the collision operator for the right hand side. This would yield the following uncoupled system
  \begin{equation*}
		\left\{ \begin{aligned}
		& \partial_t \rho  + \diverg_x (\rho \, \bm{u} ) \,=\, 0, 
		\\
		& \,
		\\
		& \partial_t(\rho \, \bm{u} ) + \diverg_x  \left(\rho \, \bm{u} \otimes \bm{u}\,+\, \rho \, T \,{\rm\bf I} \right)  \,=\, 0, 
		\\
		& \,
		\\
		& \partial_t E + \diverg_x \left(\bm{u}\, (  E \,+\, \rho \, T  ) \right) \,=\, \frac{1}{2 \,\ve} \int_{\RR^{d_v}} {\Q}\left (\mathcal M_{\rho, \bm{u}, T}, \mathcal M_{\rho, \bm{u}, T}\right )(t,x,v) \, |v|^2 \, dv.
		\end{aligned} \right.
	\end{equation*}
	One of the advantages of this uncoupled approach is that, as we will see in the following section, it generalizes the closure we used for both elastic and inelastic Boltzmann equation.
	We are currently investigating it for the case of kinetic models of flocking, such as the one presented in \cite{Fornasier:2011} (where momentum is not conserved as well).
\end{remark}
				
In the following section we give several examples of application of such a method. 

\section{Applications of the Rescaling Velocity Method}
\label{secApp}
		
We start with the simple example of the classical Boltzmann equation and then present different applications to inelastic kinetic models for granular gases.

\subsection{Application to the Boltzmann Equation}

The Boltzmann equation describes the behavior of a dilute gas of particles when the only interactions taken into account are binary	elastic collisions. It is given by equation~\eqref{eqCollision}	where $f^{\ve}(t,x,v)$ is the time-dependent particle distribution function in the phase space and the Boltzmann collision operator $\Q_\mathcal{B}$ is a quadratic operator, local in $(t,x)$
\begin{equation}
\label{eqQBoltz}
\Q_\mathcal{B} (f,f)(v) \,=\, \int_{\RR^{d_v}} \int_{\mathbb{S}^{{d_v}-1}}  B(|v-v_*|,\cos \theta) \, \left[ f'_* f' \,-\, f_* f \right] \, d\sigma \, dv_*,
\end{equation}
where we used the shorthand $f = f(v)$, $f_* = f(v_*)$, $f ^{'} = f(v')$, $f_* ^{'} = f(v_* ^{'})$.  The velocities of the colliding pairs $(v,v_*)$ and $(v',v'_*)$ are related by
\begin{equation*}
v' \,=\, \frac{v+v_*}{2} \,+\, \frac{|v-v_*|}{2} \,\sigma, \qquad v'_* \,=\, \frac{v+v^*}{2} \,-\, \frac{|v-v_*|}{2} \,\sigma.
\end{equation*}
The \emph{collision kernel} $B$ is a non-negative function which by physical arguments of invariance only depends on $|v-v_*|$ and	$\cos \theta = {\widehat u} \cdot \sigma$, where ${\widehat u} =(v-v_*)/|v-v_*|$.
			
Boltzmann's collision operator has the fundamental properties of conserving mass, momentum and energy
\begin{equation*}
\int_{\RR^{d_v}}\Q_\mathcal{B}(f,f) \, \varphi(v)\,dv \,= \,0, \qquad \varphi(v)\,=\,1,\,v,\,|v|^2
\end{equation*}
and satisfies the well-known Boltzmann's $H$-theorem (for the space homogeneous case)
\begin{equation*}
- \frac{d}{dt} \int_{\RR^{d_v}} f \, \log f \, dv \,=\, - \int_{\RR^{d_v}} \Q_\mathcal{B}(f,f) \, \log(f) \, dv \geq 0.
\end{equation*}
The functional $- \int f \log f \, dv$ is the entropy of the solution. Boltzmann $H$-theorem implies that any equilibrium distribution	function, {\it i.e.} any function which is a maximum of the entropy, has the form of a locally Maxwellian distribution \cite{CIP:94, Golse:2005}
\begin{equation} 
\label{defMaxw}
\mathcal{M}_{\rho,u,T}(v)=\frac{\rho}{(2 \, \pi \, T)^{{d_v}/2}} \exp \left( - \frac{\vert u - v \vert^2} {2 \, T} \right),
\end{equation}
where the density, velocity and	temperature of the gas $\rho$, $u$ and $T$  are computed from the distribution function as in (\ref{def:rhoE})-(\ref{def:uT}). For further details on the physical background and derivation of the Boltzmann equation we refer to~\cite{CIP:94,Villani:2002handbook}.
			
Assuming that the collision kernel is given by the \emph{variable hard spheres} model, that is,
\begin{equation*} 
B(|z|,\cos\theta) \,=\, |z|^\lambda \,b(\cos\theta)
\end{equation*}
for $\lambda \in [0,1]$ and $b$ bounded, then for a distribution $g$ given by \eqref{defRescaling}, we have using \eqref{eqQBoltz}
\begin{equation*}
	\Q_\mathcal{B}(f,f) \,=\, \omega^{\lambda - {d_v}}\, \Q_\mathcal{B}(g,g)
\end{equation*}
and the rescaled equation~\eqref{eqGcons} reads
\begin{equation}
\label{eq:g:elastic}
\frac{\partial g}{\partial t} + \diverg_x \left ( \omega \,\xi \, g\right )  - \diverg_\xi \left [ \left ( \frac{1}{\omega} \frac{\partial \omega}{\partial t} \, \xi + \xi \otimes \xi \, \nabla_x \omega \right ) g \right ] = \frac{\omega^\lambda}{\ve} \, \Q_\mathcal{B}(g,g).
\end{equation}
Moreover, applying a standard closure of the Boltzmann equation by a Maxwellian distribution~\eqref{defMaxw}, we find that the macroscopic quantities and then $\omega$ thanks to \eqref{eqWclose}, are given by
\begin{equation}
\label{eq:macro:elastic}
\left\{ \begin{array}{l}
\partial_t \rho  + \diverg_x  (\rho \, \bm{u} ) \,=\, 0, 
\\
\,
\\
\partial_t(\rho \, \bm{u} ) + \diverg_x  \left(\rho \, \bm{u} \otimes \bm{u} \,+\, \rho \, T  \,{\rm\bf I}\right) \, =\, 0, 
\\
\,
\\
\partial_t E + \diverg_x  \left(  \bm{u} \,(  E \,+\, \rho \, T )\right) \,=\, 0.
\end{array} \right.
\end{equation}

We shall now apply the rescaling method to the granular gases equation.

\subsection{Application to the Granular Gases Equation}
\label{subGranular}

A granular gas is a set of particles which interact by energy dissipating collisions leading to inelastic collisions.  This inelasticity is characterized by a collision mechanics where mass and momentum are conserved and kinetic energy is dissipated.  Thus, the collision phenomenon is a non-\emph{microreversible} process.  The velocities of the	colliding pairs $(v,v_*)$ and $(v',v'_*)$ are related by
\begin{equation*}
v'  \,=\, v \,-\, \frac{1 + e}{2}\, \left( u \cdot \sigma \right) \sigma, \qquad
v_*'\,=\, v_* \,+\, \frac{1 + e}{2} \,\left( u \cdot \sigma \right) \sigma,
\end{equation*}
where $u := v - v_*$  is the \emph{relative velocity}, $\sigma \in \SSS^{d-1}$ the \emph{impact direction} and $e \in [0,1]$ the dissipation parameter, known as \emph{restitution coefficient}. It means that the energy is dissipated in the impact direction only and will be chosen constant. The parametrization of post-collisional velocities can also be found by using some properties of the model. Indeed, it is equivalent to the conservation of impulsion and dissipation of energy:
\begin{equation} 
\label{eqMicroCons}
\left\{ \begin{array}{l}
v' \,+\, v_*' \,=\, v \,+\, v_*, 
\\
|v'|^2 \,+\, |v_*'|^2 \,-\, |v|^2 \,-\, |v_*|^2  \,=\, - \displaystyle\frac{1-e^2}{2} \,| u \cdot \sigma |^2 \,\leq\, 0.
\end{array} \right.
\end{equation}
This microscopic mechanism allows us to describe the granular collision operator $\Q_\mathcal{I}$. If the nonnegative particles distribution function $f$  depends only on $v \in \mathbb{R}^{d_v}$, the collision operator $\Q_\mathcal{I}(f,g)$ can be expressed in the following weak form: given a smooth test function $\psi$,
\begin{equation} 
\label{Qweak}
\int_{\RR^{d_v}} \Q_\mathcal{I}(f,f)\, \psi \, dv \,:=\, \frac{1}{2 }\int_{\mathbb{R}^{d_v} \times \mathbb{R}^{d_v} \times \mathbb{S}^{{d_v}-1}} |v-v_*| \,f_{*} \, f \, \left(\psi' \,+\, \psi_*' \,-\, \psi \,-\, \psi_* \right)\, b(\cos \theta) \,d\sigma \, dv \, dv_*,
\end{equation}
where we have used the shorthand $\psi' := \psi(v')$, $\psi_*' := \psi(v_*')$, $\psi := \psi(v)$ and $\psi_* := \psi(v_*)$. 

Taking $\psi(v) = 1$, $v$ and $|v|^2$ in \eqref{Qweak}, the relations \eqref{eqMicroCons} yield conservation of mass and momentum at the kinetic level and  dissipation of kinetic energy
\begin{equation*}
\int_{\RR^{d_v}} \Q_\mathcal{I}(f,f)\, |v|^2 \,dv \, =\, -\,D(f).
\end{equation*}
where $D(f)$ is the \emph{energy dissipation functional}, given by
\begin{equation}
D(f)  \,=\, C\, (1-e^2)\, \iint_{\mathbb{R}^{d_v} \times \mathbb{R}^{d_v}}f\, f_* \,|v-v_*|^3 dv \, dv_* \,\geq\, 0, 
\label{eqDissipQ}
\end{equation}
for $C$ an explicit, nonnegative constant depending on the cross section $b$ and the dimension.  The decrease of the energy, together with conservation of mass and momentum imply that the equilibria of the collision operator are Dirac distributions $\rho \, \delta_{\bm u}(v)$  where the density $\rho$ and momentum $\rho \bm u$ are prescribed. This is the major problem for the numerical study of equation~\eqref{eqCollision} with deterministic algorithms because of this concentration in velocity variables. In that case, the rescaling method (\ref{defRescaling}) becomes an essential tool to follow the concentration in velocity of the distribution function.

For the hard sphere kernel (\ref{Qweak}) and for a rescaled distribution function $g$ given by~\eqref{defRescaling}, we have
\begin{equation*}
\Q_\mathcal{I}(f,f) = \omega^{1-{d_v}} \, \Q_\mathcal{I}(g,g).
\end{equation*}
The rescaled equation~\eqref{eqGcons} then reads
\begin{equation} 
\label{eqGgran}
\frac{\partial g}{\partial t} \,+\, \diverg_x \left ( \omega \,\xi \, g\right )  \,-\, \diverg_\xi \left [ \left ( \frac{1}{\omega} \frac{\partial \omega}{\partial t} \, \xi + \xi \otimes \xi \, \nabla_x \omega \right ) g \right ] \,=\, \frac{\omega}{\ve} \, \Q_\mathcal{I}(g,g).
\end{equation}
To complete the rescaling procedure, it remains to derive a set of equations for $\rho$ and $E$ to define the rescaling function $\omega$. To this end, we will assume that the restitution coefficient $e$ is close to $1$, which corresponds to the weak inelasticity case, in order to write the collision operator as a perturbation of the classical Boltzmann operator and then write the system~\eqref{eqSysConsUncoupled} associated to granular gases.  Let us mention that there are alternative methods for deriving hydrodynamic equations from dissipative kinetic models, which  can be found in the review paper of M. Bisi, J.A. Carrillo and G. Spiga \cite{Bisi:2010}. Here, we focus on the approach developed by G. Toscani \cite{Toscani:2004, Toscani:2007}, which is reminiscent of a the classical works of Brilliantov, Goldhirsch, Jenkins, P\"oschel and Zanetti \cite{brilliantov:2004, Goldhirsch:1993, Jenkins:1985}:  we write the granular gases operator as a perturbation of the elastic Boltzmann operator
\begin{equation}
\Q_\mathcal{I}(f,f) \,=\,  \Q_\mathcal{B}(f,f) \,+\, \ve \, \mathcal I(f,f) \,+\, \mathcal O(\ve^2), 
\label{QPerturb}
\end{equation}
where $\mathcal I$ is a quadratic, energy dissipative nonlinear friction operator given by
\begin{equation*}
  \mathcal I(f,f)(v) \,=\, \frac{1}{2} \diverg_v \left ( \int_{\mathbb{R}^{d_v} \times \mathbb{S}^{{d_v}-1}}  |u \cdot \sigma| \, (u \cdot \sigma) \, \sigma \, f_{*} \, f\, b(\cos \theta) \, d\sigma \, dv_*\right ).
\end{equation*}
		
Using the expansion \eqref{QPerturb} of the operator with respect to $\ve$, it was shown in \cite{Toscani:2004} that we can derive the macroscopic system~\eqref{eqSysConsUncoupled} associated to $\Q_\mathcal{I}$. Assuming that the inelasticity is of the same order of the Knudsen number:
\begin{equation} 
\label{quasi}
1-e \,=\, \varepsilon,
\end{equation}
then, up to the second order in $\ve$, the macroscopic quantities are solution to
\begin{equation} 
\label{eqSysCons}
\left\{ \begin{array}{l}
\partial_t \rho \,+\, \diverg_x  (\rho \, \bm{u}) \,=\, 0, 
\\
\,
\\
\partial_t(\rho \, \bm{u}) \,+\, \diverg_x  \left(\rho \,( \bm{u} \otimes \bm{u})  \,+\, \rho \, T \,{\rm\bf I}\right) \,=\, 0, 
\\
\,
\\
\partial_t E \,+\, \diverg_x  \left( \bm{u}  \,( E\,+\, \rho \, T ) \right) \,=\, -K_{{d_v}} \, \rho^2 \, T^{3/2},
\end{array} \right.
\end{equation}
for an explicit and nonnegative constant $K_d$ \cite{Toscani:2004}. The solution of this system will give the scaling function $\omega$ in~\eqref{eqWclose}.
			
%%%%% FF %%%
\begin{remark}
Let us notice that during the evolution if hydrodynamic quantities $(\rho,u,T)$ given by the solution to (\ref{eqSysCons}) and the ones obtained by computing the moments of the distribution function $g$ solution to (\ref{eqGgran}) are completely different, it is possible to update the scaling, that is the hydrodynamic model, starting with new macroscopic quantities which corresponds to moments of the distribution function $g$. If this happens it only means that the velocity grid in not well fitted with the scale of the distribution function $g$ and this process corresponds to adapt the grid in the velocity space. 
\end{remark}
		
\subsection{Application to a Simplified Granular Gases Model}
\label{subBGK}

The rescaling velocity method can be also applied to a simplified model designed to mimic the behavior of the granular gases equation with a simpler collision operator. Indeed, A. Astillero and A. Santos derived in~\cite{Astillero:2004,Astillero:2005} a very simple approximation of the granular gases equation, containing its main features. It consists in replacing the \emph{bilinear} collision operator $\Q_\mathcal{I}$ by a \emph{nonlinear} relaxation operator, supplemented with a drag force in the velocity space, depending on the moment of the particle distribution function itself, in order to match the correct macroscopic dissipation of kinetic energy of the granular gases equation. 

The relaxation operator we will use is the so-called Ellipsoidal Statistical BGK operator (ES-BGK) \cite{AndriesPerthame:2000}. To this aim we  first define some macroscopic quantities of the particle distribution function $f$ such as the opposite of the \emph{stress tensor}
\begin{equation*}
\Theta_f(t,x)\, =\, \frac{1}{\rho_f} \int_{{\RR}^{d_v}}  (v-\bm{u}_f)\otimes (v-\bm{u}_f)\,f(t,x,v) \,dv.
\end{equation*} 
Therefore the \emph{translational temperature} is related to the stress tensor as $T_f = {\rm tr}(\Theta_f)/d_v$. We finally introduce the corrected tensor
\begin{equation*}
\Tau_f(t,x) \,\,=\,\, \left[(1-\nu) \, T_f \,{\rm\bf I} \,\,+\,\,\nu \,\Theta_f\right](t,x),  
\end{equation*} 
which can be viewed as a linear combination of the initial stress tensor $\Theta_f$ and of the isotropic stress tensor $T_f \,{\rm \bf I}$ developed by a Maxwellian distribution. The parameter $-\infty< \nu < 1$ is used to modify the value of the Prandtl number through the formula 
$$
0 \,\leq\, {\rm Pr}  \,=\, \frac{1}{1-\nu} \,\leq\, +\infty\quad {\rm for } \quad\nu\in (-\infty\,,\,1).
$$
Thus we introduce a corrected Gaussian $\G[f]$ defined by
$$
\G[f]= \frac{\rho_f}{\sqrt{{\rm det}(2\pi\,\Tau_f)}}\,\exp\left(-\frac{(v-\bm{u}_f)\,\Tau_f^{-1}\,(v-\bm{u}_f)}{2}\right)
$$ 
and the corresponding collision operator is now
\begin{equation} 
\label{eqOpFPBGK}
\Q_\mathcal{S}(f) =  \tau(P_f) \,\left( \G[f] \,-\, f \right) \,\,+\,\, C \,\rho_f \,T_f^{1/2} \, \diverg_{v} \left ( (v-\bm{u}_f) f\right ),
\end{equation}
where $\tau$ is a given function depending on the \emph{kinetic pressure} $P_f:= \rho_f \, T_f $ and $C$ is given by 
$$
C\,=\,(1-e)\, \frac{K_{d_v}}{d_v}.
$$  
The coefficient $e$ corresponds to the restitution coefficient of the granular operator $\Q_\mathcal{I}$ and will be assumed equal to $1-\ve$ as in~\eqref{quasi}.  This expression yields the correct granular energy dissipation \eqref{eqDissipQ}.

This new operator replaces the inelastic collision operator by an ``equivalent'' elastic operator and the drift term acts like an external drag force which mimics the cooling effect of the particles. In particular, the hydrodynamic equations associated to the rescaled equation~\eqref{eqGcons} with operator~\eqref{eqOpFPBGK} are exactly given by the system~\eqref{eqSysCons}. The operator $\widetilde{\Q}_\mathcal{S}$ such that $\widetilde{\Q}_\mathcal{S}(g) = \Q_\mathcal{S}(f)$ is given by
\begin{equation*}
\widetilde{\Q}_\mathcal{S}(g) \,=\, \frac{1}{\omega^{d_v}}\left [ \tau_g \,\left(\G[g] - g\right) \,\,+\, \,C \,\omega \,\rho_g \,T_g^{1/2} \, \diverg_{\xi} \left (\left (\xi\,-\, \bm{u}_g \right )\, g\right )\right ],
\end{equation*}
using the shorthand $\tau_g := \tau\left ( \omega^2 \, P_g \right )$. Equation~\eqref{eqGcons} then simply reads
\begin{equation*} 
%\label{eqGFP}
\frac{\partial g}{\partial t} \,+\, \diverg_x  \left ( \omega \,\xi \, g\right ) \,=\, \frac{ \tau_g }{\ve} \, (  \G[g] \,-\, g )\,+\, \diverg_{\xi} \left [ \left ( C\, \omega\, \rho_g \,T_g^{1/2}\, \left (\xi\,-\, \bm{u}_g \right ) \,\,+\,\, \frac{1}{\omega} \frac{\partial \omega}{\partial t} \, \xi + \xi \otimes \xi \, \nabla_x  \omega \right ) g \right ] ,
\end{equation*}
where the function $\omega$ is given by~\eqref{eqWclose} using the closed system~\eqref{eqSysCons}.

\section{Boundary Conditions} 
\label{secBoundCond}	
\begin{figure}					
\begin{center}
\includegraphics{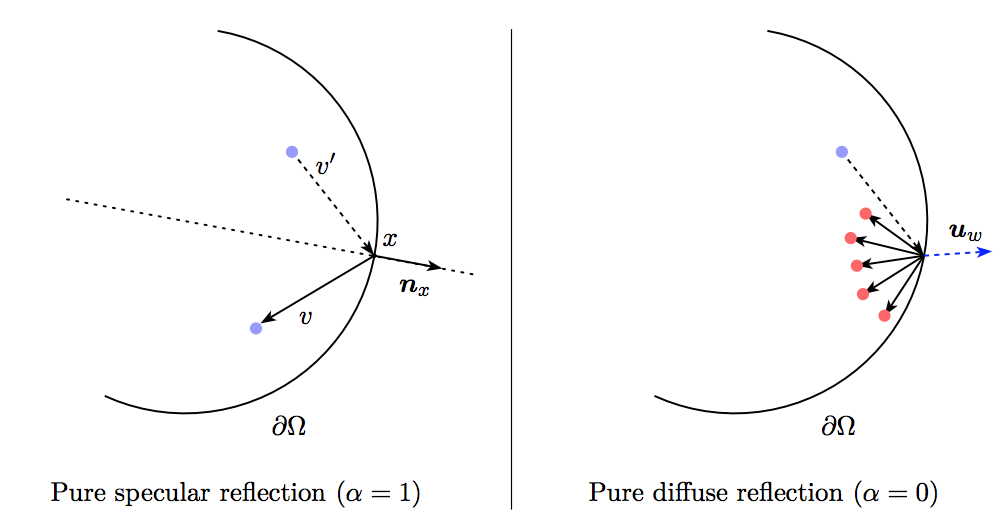}
\caption{Schematic representation of Maxwell conditions at the boundary.\label{figMaxwBound}}
\end{center}
\end{figure}

In order to define completely the mathematical problem for the kinetic equation \eqref{eqCollision}, one has to set boundary conditions on $\partial \Omega$. The most simple (yet physically relevant) are the so-called Maxwell boundary conditions~\cite{Maxwell:1867}, where a fraction $1-\alpha \in [0,1]$ of the particles is \emph{specularly} (elastically) reflected by a wall, whereas the remaining fraction $\alpha$ is \emph{thermalized}, leaving the wall (moving at a velocity $\bm u_w$) at a temperature $T_{w} > 0$ (see Figure~\ref{figMaxwBound}).
The parameter $\alpha$ is called \emph{accommodation} coefficient, because it expresses the tendency of the gas to accommodate to the state of the wall. More details on the subject can be found in the textbook of Cercignani~\cite{Cercignani:1988}.
 The main difficulty consists to apply appropriate boundary conditions on the macroscopic model which are consistent with the kinetic boundary value problem.
 
Let us consider a distribution function $f(t,x,v)$ solution to \eqref{eqCollision} for $(t,x,v) \in \RR_+ \times \Omega \times \RR^{d_v}$ with $\Omega$ a  bounded Lipschitz-continuous domain of $\RR^{d_x}$. We assume the  boundary $\partial \Omega$ to have a unit outer normal $\bm{n}_x$ for all $x \in \partial \Omega$.  Then, for an \emph{incoming} velocity $v \in \RR^{d_v}$ (namely $v \cdot \bm{n}_x \leq 0$), we set
\begin{equation*}
f(t,x,v) \,=\, \alpha \,\mathcal{R}f(t,x,v) \,+\, (1-\alpha) \,\mathcal{M}f(t,x,v), \quad \forall x \in \partial \Omega, \, v\cdot \bm{n}_x \leq 0,
\end{equation*}
with
\begin{equation} 
\label{eqBoundCond}
\left\{ \begin{array}{l}
\mathcal{R}f(t,x,v) \,\,=\,\, f(t,x, v\,-\, 2 (\bm{n}_x \cdot v)\, \bm{n}_x), 
\\
\,
\\
\mathcal{M}f(t,x,v) \,\,=\,\, \mu(t,x) \,\mathcal {M}_{1, \bm{u}_w, T_w}(v).
\end{array}\right.
\end{equation}
In this last expression, $\mathcal {M}_{1,  \, \bm{u}_w,  \, T_w}$ is the wall Maxwellian
\begin{equation*}
\mathcal {M}_{1,  \, \bm{u}_w, \,  T_w}(v) := \frac{1}{(2 \pi \, T_w)^{d_v/2}} \exp\left (-\frac{|v - \bm u_w|^2}{2\,T_{w}}\right )
\end{equation*}
and $\mu$ insures global mass conservation:
\begin{equation*} 
\int_{v\cdot \bm{n}_x \geq 0} f(t,x,v) \,v \cdot \bm{n}_x \, dv \,\,=\,\, - \mu(t,x) \int_{v\cdot \bm{n}_x < 0} \mathcal {M}_{1,  \, \bm{u}_w,  \, T_w}(v) v \cdot \bm{n}_x \, dv.
\end{equation*}
Thanks to these definitions, it is classical (see for example \cite{Cercignani:1988}) that the following theorem holds:

\begin{theorem}\label{thmConsBound}
  Let $f$ be a smooth solution to the kinetic equation \eqref{eqCollision} with Maxwell boundary conditions \eqref{eqBoundCond}. Then,
  \begin{enumerate}
    \item the mass of $f$ is globally conserved for any accommodation coefficient $\alpha \in [0,1]$;
    \item if $\Q$ preserves the kinetic energy, then in the specular case $\alpha = 1$ the energy of $f$ is globally conserved ;
    \item if $\Q$ produces entropy, namely \[ -\int_{\RR^{d_v}} \Q_\mathcal{B}(f,f) \, \log(f) \, dv \geq 0\] then the entropy  functional $\int f \log f \, dv \, dx$ is dissipated for any accommodation coefficient $\alpha \in [0,1]$.
  \end{enumerate}
\end{theorem}
%
%\begin{proof}
%  (1) We start by multiplying equation \eqref{eqCollision} by $\varphi$ and integrating in velocity and space variables. We have made the general assumption that $\Q$ preserves the mass. Then, using the Stokes formula   and taking  $\Psi (f) = f$ and $\varphi(v) = 1$ (which is symmetric with respect to the normal component of $v$) in \eqref{eqSpecularInt} is enough to guarantee the mass conservation for the specular boundary conditions. The relation \eqref{eqDiffusInt} allows to conclude in the diffusive case.
%  
%  (2) In the pure specular case, the elastic collisions at the boundary yield the microscopic conservation of kinetic energy: $|v'|^2 = |v|^2$. Then, taking  $\Psi (f) = f$ and $\varphi(v) = |v|^2$ (which is symmetric with respect to the normal component of $v$) in \eqref{eqSpecularInt} and using the Stokes formula gives the conservation of kinetic energy.
%  
%  (3) Multiplying equation \eqref{eqCollision} by $\log f + 1$ and taking  $\Psi (f) = f \, \log f$, $\varphi(v) = 1$ in \eqref{eqSpecularInt} yields the dissipation of entropy, thanks to Stokes formula.
%\end{proof}

We shall present the boundary conditions associated to the Maxwell conditions for the rescaled distribution $g$ and the hydrodynamics quantites $(\rho, \rho \, \bm u, E)$ solution to the  Euler' system. We split in two parts, pure specular and pure diffusive boundary conditions.

\subsection{Specular Boundary Conditions}

The natural translation of a specular reflection of a particle of velocity $v$ in the physical space, at a point $x$ of the boundary $\partial \Omega$ is a specular reflection at the point $x$, but with a rescaled velocity $\xi = v \, \omega$. More precisely, for an incoming velocity $\xi$ (namely $\xi \cdot \bm{n}_x \leq 0$), we set 
\begin{equation*}
{\mathcal{R}} \, g(t,x,\xi) =  g\left (t,x, \left [\xi- 2 (\bm{n}_x \cdot \xi) \, \bm{n}_x\right ]\right ).
\end{equation*}

Concerning the macroscopic quantities, we need to impose at the fluid level the so-called \emph{slip boundary} condition \cite{sone:2007}. It corresponds to a zero-flux boundary conditions for the mean velocity:
\[ \bm u \cdot \bm{n}_x = 0. \]
This condition yields in particular the global conservation of mass and kinetic energy that we observed at the kinetic level in Corollary \ref{thmConsBound}. It then remains to fix a condition for the momentum flux. To mimic the specular reflection, we reverse the velocity at the boundary, as presented in \cite{Besse:2011}.

\subsection{Diffusive Boundary Conditions}	

We want to impose a velocity $\bm u_w$ and a temperature $T_w$ at a point $x$ of the boundary, using a wall maxwellian $\mathcal {M}_{1,  \, \bm{u}_w,  \, T_w}$ in the physical space. Then, for an incoming velocity $\xi$ in the rescaled space, we set
\begin{equation*}
{\mathcal{M}} \, g(t,x,\xi) \,\,=\,\, {\mu}(t,x) \,\mathcal{M}_{1, \, \widetilde{\bm{u}}_w, \, \widetilde{T}_w}\left (\xi\right ).
\end{equation*}
In this last expression, $\mathcal{M}_{1, \, \widetilde{\bm{u}}_w, \, \widetilde{T}_w}$ is the ``rescaled'' wall Maxwellian, with $\widetilde{\bm{u}}_w = \omega \, \bm{u}_w$, $\widetilde{T}_w := \omega^2 \, T_w$ and ${\mu}$ insures the global mass conservation:
\begin{equation*}
\int_{\xi \cdot \bm{n}_x > 0} g(t,x,\xi) \,\xi \cdot \bm{n}_x \, d\xi \,\,=\,\, - \mu(t,x)\, \int_{\xi \cdot \bm{n}_x < 0} \mathcal{M}_{1, \,  \widetilde{\bm{u}}_w, \, \widetilde{T}_w}\left ( \xi\right )\, \xi \cdot \bm{n}_x \, d\xi.
\end{equation*}

At the macroscopic level, we impose a zero-flux boundary conditions for the mean velocity:
\begin{equation*}
\bm{u} \cdot \bm{n}_x = 0.
\end{equation*}
This yields the global conservation of mass for the macroscopic model, corresponding to the result of Corollary \ref{thmConsBound}. We also want to fix the temperature $T_{w} >0$ of the wall. We set
\begin{equation*}
E(x) \,=\, \frac{{d_v}}{2}\, \rho(x)\, T_{w}.
\end{equation*}
These boundary conditions allow to solve the boundary value problems  (\ref{eq:macro:elastic}) or  (\ref{eqSysCons}).		

\begin{remark}
Of course, the boundary conditions we  presented here for the macroscopic model are very approximate, and do no reflect exactly what type of conditions have to be imposed for this model to be the true hydrodynamic limit of the kinetic equation \eqref{eqCollision}. 
An analysis of the boundary layer has to be carried out properly in order to do so, as was done for example by F. Coron in \cite{Coron:89} for the elastic Boltzmann equation.
Nevertheless, we will see in the numerical experiments that our conditions are  accurate enough for the rescaling velocity method to behave properly, because we don't need the exact hydrodynamic limit of our problem to compute the scaling function.
\end{remark}

\section{Numerical Algorithms}
\label{secAlgo}
\subsection{Discretization of the Problem}
\label{subDiscret}

We will briefly present  the numerical schemes we used for the discretization of the uncoupled models \eqref{eqGcons} and \eqref{eqSysConsUncoupled}.
We need to distinguish three independent problems: the treatment of the free transport term, the collision term and the system of conservation laws.

\vspace{0.2cm}

\paragraph{\textbf{The finite volume method for a conservative transport equation.}}
  
  The first numerical problem concerns the space discretization of the conservative transport equation
  \begin{equation}
	  \label{eqFreeTransp}
	  \left \{ \begin{aligned}
	   & \frac{\partial g}{\partial t} + \diverg_X \left ( a(t,X) \, g\right ) = 0, \ \forall \, (t,X) \in \RR_+ \times \mathcal O, \\
	   & \, \\
	   & g(0,X) = g_0(X), 
	   \end{aligned} \right.
	\end{equation}
	for a smooth function $a : \RR_+ \times \mathcal O \to \RR^{d_X}$ and a Lipschitz-continuous domain $\mathcal O \subset \RR^{d_X}$. 
	Our treatment of the problem is made in the framework of finite volume schemes, on Cartesian grid. This type of methods is particularly well suited to solve conservative transport equation. We choose to use an upwind method with the so-called Van Leer's flux limiter, as presented in the classical paper \cite{VanLeer:1977a}.

\vspace{0.2cm}
	
\paragraph{\textbf{The spectral method for the collision operator}.}
  
  We now focus on the numerical resolution of
  \begin{equation} \label{eqCollHomog}
    \frac{\partial g}{\partial t} = \Q(g , g ),
  \end{equation}   
  where $\Q = \Q_\mathcal{B}$ or $\Q_\mathcal{I}$. 
  
  Any deterministic numerical method for the computation of the Boltzmann operator requires to work on a \emph{bounded} velocity space. Our approach consists in adding some non physical binary collisions by periodizing the particle distribution function and the collision operator, and then using Fourier series to compute the truncated operator. 
  This implies the loss of some local invariants, but a careful periodization allows at least the preservation of mass. This periodization is the basis of spectral methods, that we will use in our numerical simulations.

	Fourier techniques for the resolution of the Boltzmann equation have been first introduced independently by L. Pareschi and B. Perthame in \cite{pareschi:1996} and by A. Bobylev and S. Rjasanow in \cite{Bobylev:1997}. 
	It has since been investigated by a lot of authors. One can find numerous results about this method in the article \cite{PaRu:SINUM:2000} of L. Pareschi and G. Russo. It has been recently shown to be convergent by F. Filbet and C. Mouhot in \cite{filbet:2011Conv} for the elastic case and allows spectral accuracy \cite{canuto:88}. 
	It has then been derived for the granular gases operator \eqref{Qweak} by F. Filbet, L. Pareschi and G. Toscani in \cite{filbet:2005}. 
  
	 The alternative approach of \cite{Bobylev:1997} has been developed in $3d_v$ and improved in \cite{Gamba:2009,gamba:2010}.
	 The concept here is ``dual'' from the previous approach, since it consists in first using Fourier transform of the collision operator and then truncating the Fourier transformed operator. 
	 Another difference for this setting is that conservation of momentum and kinetic energy are achieved by solving a constrained optimization problem, which is not the case of the method we use in the following of this article.
\vspace{0.2cm}

\paragraph{\textbf{The WENO approach for a system of conservation laws.}} 

	The fluid equations \eqref{eqSysCons} can be written as a system of nonlinear hyperbolic conservation laws with source term
	\begin{equation*}
	  \frac{\partial U}{\partial t} + \diverg_x F(U) = G(U),
	\end{equation*}
	where $U:\mathbb{R}_+ \times \Omega \to \mathbb{R}^n$, $F,G:\mathbb{R}^n \to \mathbb{R}^n$ and $n={d_v}+d_x$. There exists a large number of numerical methods to solve this type of problem, and we choose to use finite differences. 
	
	The reconstruction of the fluxes will be made by a high-order Weighted Essentially Non-Oscillatory (\emph{WENO}) method with a $5$ points stencil. Indeed, it is very well adapted to the treatment of shocks in fluid systems. The reader can consult the lecture notes \cite{shu:1998} of C.W. Shu for most of the details about this method.

  The rescaling function $\omega$ (and its time and space derivatives) given by the macroscopic scaling~\eqref{eqWclose} is then evolved in time thanks to the WENO fluxes $\mathcal{F}_{i+\frac 12}$.

We shall now present numerical simulations of the rescaling velocity method for energy dissipative equations.

\section{Numerical Simulations} 
\label{secNumRes}

\subsection{Test 1: Convergence toward a Dirac mass}
In this section, we will deal with a space homogeneous distribution $f=f(t,v)$, solution to the granular gases equation					
\begin{equation} 
\label{eqboltzHomog}
\left\{
\begin{array}{l}
\displaystyle \frac{\partial f}{\partial t}  = {\Q_\mathcal{I}(f, f)},
\\
\,
\\
f(t=0) = f_0,
\end{array}\right.
\end{equation}
with an initial condition chosen as the sum of two gaussian distributions of total temperature equal to $1$:
\begin{equation*}
f_{0}(v) := \frac{1}{2 \, ( 2 \, \pi \, \sigma^2)^{{d_v}/2}}\left ( \exp\left (-\frac{|v - \sigma \, \bm e |^2}{2 \, \sigma^2}\right ) + \exp\left (-\frac{|v + \sigma \, \bm e|^2}{2 \, \sigma^2}\right ) \right ), \quad \forall \, v \in \RR^{d_v},
\end{equation*}
where $\bm e = (1, \ldots, 1)^\mathsf{T}$ and $\sigma = 1/10$. With this choice, we have
\begin{equation}
\label{propCIhomog}
  \int_{\RR^{d_v}} f_0(v)\, 
  \left( \begin{array}{c} 1 \\ v \\  |v|^2 \end{array}\right) dv \,\,=\,\, \left( \begin{array}{l} 1 \\ 0\\d_v\end{array}\right).  
\end{equation}
This initial datum is far from equilibrium, for all restitution coefficients $e \in [0,1]$.
The corresponding scaled distribution $g=g(t,\xi)$ given by (\ref{defRescaling}) is then solution to the homogeneous granular equation with a time-dependent drift term
\begin{equation} 
\label{eqGhomog}
\left \{ \begin{aligned}
  & \frac{\partial g}{\partial t} - \nabla_\xi \cdot \left ( \frac{1}{\omega} \frac{d \omega}{dt} \, \xi \, g \right ) \,=\, \omega\, \, \Q_\mathcal{I}(g,g),
	\\
	& \,
	\\
	& g(t=0) = g_0.
\end{aligned}\right.
\end{equation}
The initial condition $g_{0}$ is chosen by $g_0=f_{0}$ because $\omega(0) = 1$ according to \eqref{propCIhomog}.

To study the evolution of $\omega$, we propose to compare the rescaling approach proposed by F. Filbet and G. Russo in \cite{Filbet:2004}  and the macroscopic one based on the closure of the kinetic model. In the space homogeneous setting, the mass $\rho = 1$ is preserved. Hence, we have \[ \omega(t) = \sqrt{\frac{2}{{d_v}} E(t)}\]
and it is enough to compute the time evolution of the kinetic energy $E$ to obtain $\omega$.
Then on the one hand, following  \cite{Filbet:2004}, we will study the distribution $g$ solution to \eqref{eqGhomog} coupled with the exact kinetic energy
\begin{equation} 
\label{eqEexact}
\frac{d E}{dt} \,=\, \omega^3 \, \int_{\mathbb{R}^{d_v}} {\Q_\mathcal{I}(g, g)} \,\frac{|\xi|^2}{2} \, d\xi. 
\end{equation}
On the other hand, from the result of Subsection \ref{subGranular}, we replace the ordinary differential equation \eqref{eqEexact} by 
\begin{equation} 
\label{eqEapprox}
\frac{dE}{dt} \,=\, -K_{{d_v}}\, \rho^2 \,T^{3/2},
\end{equation}
where ${d_v} \,T(t) \,=\, 2\, E(t)$ because of the conservation of momentum.  Let us emphasize that \eqref{eqGhomog}-\eqref{eqEapprox} is an  uncoupled model. We will then compare it with the solution to \eqref{eqboltzHomog} in classical variables.

\vspace{0.2cm}

\paragraph{\bf Description of the results.}

We start with a space homogeneous test problem with $d_v=2$, designed to compare the results obtained with both coupled \eqref{eqGhomog}-\eqref{eqEexact} and uncoupled models \eqref{eqGhomog}-\eqref{eqEapprox}, in the mildly inelastic regime $e=0.8$. 
With this test, we intend to show that, even if we are not computing exactly the kinetic energy (provided that we are not in a quasi-elastic setting $e \sim 1$), the uncoupled model will yield exactly the same results than the coupled one even for an initial datum far from equilibrium.
We choose the computational domain in velocity as the box $\mathcal{D}_V=[-V,V]^{d_v}$ with $V = 9$ to avoid some possible loss of mass due to the drift term. We take successively $N_v = 32$ and $N_v = 48$ half Fourier modes in each direction of the velocity space.

\begin{figure}
\begin{center}
\includegraphics{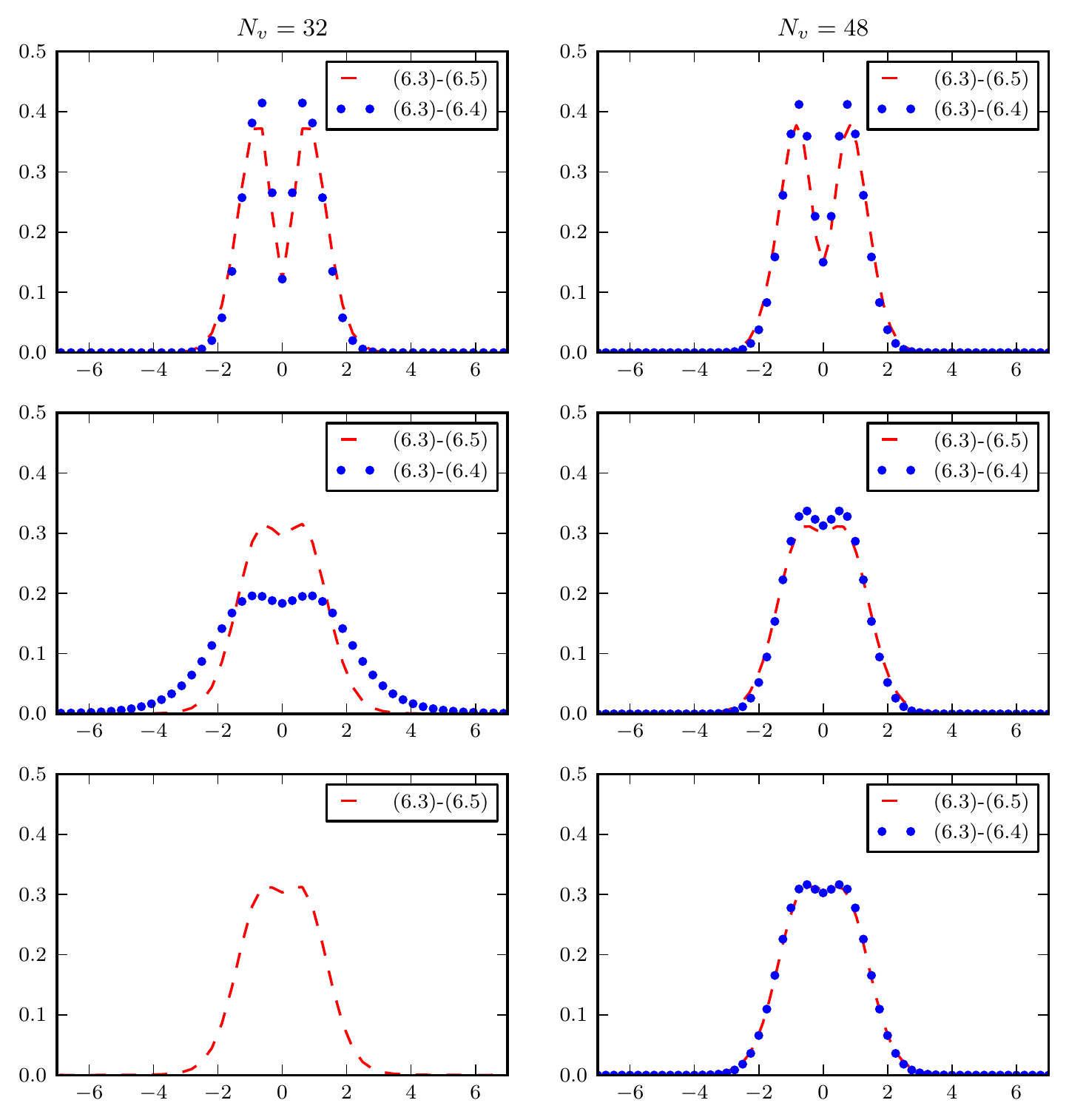}
\caption{Test 1 -- Comparison of the first marginal distribution $g^{(1)}$ for coupled model \eqref{eqGhomog}-\eqref{eqEexact} (points) and uncoupled model \eqref{eqGhomog}-\eqref{eqEexact} (dashed line), with $e = 0.8$, at times $t= 0.5$, $10$ and $50$ for $N_v = 32$ (left) and $N_v = 48$ (right). }
\label{figHomogCompG}
\end{center}
\end{figure}

We present in Figure \ref{figHomogCompG}  the time evolution of the first marginal distribution
\[ g^{(1)}(t, \xi_1) := \int_{\RR} g(t, \xi_1, \xi_2) \, d\xi_2. \]
We can see that both models give close results for short time in the coarse and fine grids, but some oscillations eventually appear in the coupled model, leading to the breakout of the scheme. 
This is due to the wrong computation of the kinetic energy of $g$ when we take a sparse grid in the velocity space, because of the heavy tails of this distribution \cite{Mischler:20062}.
This phenomenon does not occur anymore when we refine the grid. 
In this case, the results given by both models are in very good agreement, especially for the long time behavior: the equilibrium distribution is very well captured by the new uncoupled model.

\vspace{0.2cm}
\noindent

We then compare the results of this test for the uncoupled model with a solution to equation \eqref{eqboltzHomog} computed in classical variables with the spectral scheme.
We take $V=8$ and $N_v = 32$ half modes in each directions of the classical velocity space. 

\begin{figure}
	\begin{center}
		\includegraphics{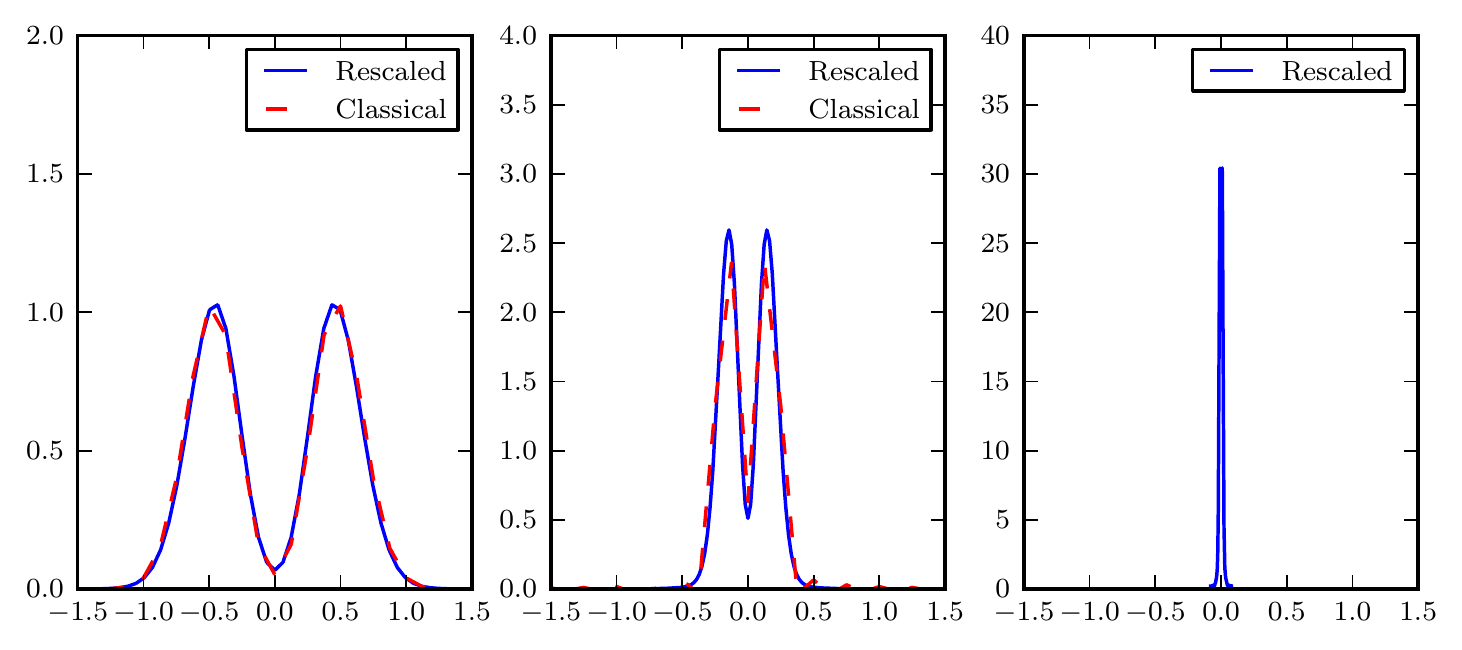}
		\caption{Test 1 -- Comparison of the first marginal $f^{(1)}$ of the solution to \eqref{eqboltzHomog} computed in classical (dashed line) and rescaled variables (solid line), with $e=0.8$, at times $t=0.3$, $1$ and $20$.}
		\label{figHomogF}
	\end{center}
\end{figure}

We observe in Figure \ref{figHomogF} that both models give the same result for short time. Then, on the one hand, the solution obtained by the spectral scheme in classical variables breaks down for intermediate time ($t \sim 1$) because of the concentration which yields spurious oscillations. We can see in particular the loss of symmetry of the solution at time $t=1$. On the other hand, the rescaled variables prevent concentration and thus allows the spectral scheme to achieve its full accuracy: the Dirac mass is very well captured for large time. This is also due to the correct change of scales in the support of the scaled distribution, thanks to the expression of $\omega$. 

In this numerical example, we can also notice that the solution to the space homogeneous equation exhibits a clear self-similar behavior: it converges more quickly towards a self-similar solution made of two ``stretched'' bumps than towards the equilibrium solution of the equation, which is a Dirac mass. This observation is perfect agreement with the theoretical results of \cite{mischler:20091}. This is even true when we choose an initial condition far from equilibrium and a mild restitution coefficients, results which are up to our knowledge not proved yet. It is important to insist on the robustness of this rescaling technique, which is simply based on the temperature dissipation estimate \eqref{eqEapprox}, since the initial datum is far from a Maxwellian distribution and the restitution coefficient ($e\simeq 0.8$) is not close to one.

\subsection{Test 2: Trend to Equilibrium in the Non Homogeneous Case}

We consider now the full inelastic Boltzmann equation in dimension ${d_v}$ on the torus
\[
\frac{\partial f}{\partial t} + v \cdot \nabla_x f = \frac{1}{\ve}\,\Q_{\mathcal{I}}(f,f), \quad x\in \TT^{d_x}, v \in \RR^{d_v},
\]
in the physical case of hard spheres collision kernel with a constant cross section $b$, and successively periodic and specular boundary conditions in $x$.  We first introduce the hydrodynamical fields associated to a kinetic distribution $f(t,x,v)$. These are the $({d_v}+ 1)$ fields of density $\rho$ (scalar) and mean velocity $\bm u$ (vector valued)  defined by the formulas (\ref{def:rhoE})-(\ref{def:uT}). Whenever $f(t,x,v)$ is a smooth solution to the granular gases equation with periodic boundary conditions, one has the global conservation laws
for mass and momentum 
\begin{eqnarray*}
&& \frac{d}{dt}\int_{\TT^{d_x}\times\RR^{d_v}} f(t,x,v)\,dx \, dv = 0, \\
&& \frac{d}{dt}\int_{\TT^{d_x}\times\RR^{d_v}} f(t,x,v)\,v\,dx \, dv = 0.
\end{eqnarray*}
Therefore, without loss of generality we shall impose 
\begin{eqnarray*}
\int_{\TT^{d_x}\times\RR^{d_v}} f(t,x,v) \,dx \, dv\,\,=\,\,1, \quad\int_{\TT^{d_x}\times\RR^{d_v}} f(t,x,v)\,v \,dx \, dv\,\,=\,\,0,
\end{eqnarray*}
where we normalized the measure on the torus: $\textrm{m}\left (\TT^{d_x}\right ) = 1$. These conservation laws are then enough to uniquely determine the stationary state of the inelastic Boltzmann equation: the normalized global Dirac distribution in velocity and constant in space
\begin{equation}
\label{maxwglob}
f_\infty(v) = \delta_0\left(v\right).
\end{equation}
Our goal here is to investigate numerically the long time behavior of the solution $f$. If $f$ is any reasonable solution of the inelastic Boltzmann equation, satisfying certain {\em a priori} bounds of compactness, then it is expected that $f$ does  converge to the global distribution function $f_\infty$ as $t$ goes to $+\infty$. The point is that on the opposite to the elastic Boltzmann equation, the solution does not satisfy any entropy theorem. Therefore, to measure the convergence to equilibrium we compute first the evolution of the global kinetic energy, then the evolution of the deviation from the global mass 
\[
\mathcal L(t) = \|\rho(t, \cdot) - 1\|_{L^1}.
\] 

We present numerical results for the full granular gases equation \eqref{eqCollision}, using the uncoupled rescaled model \eqref{eqGgran}-\eqref{eqSysCons}. The initial condition $g_0(x, \xi) = \omega(0, x)^{d_v} f_0\left (x, \omega(0, x) \, \xi\right )$ is set as
\begin{equation} \label{eqCItrendToEqu}
  \left \{ \begin{aligned}
    & f_{0}(x,v) = \frac{\rho_{0}(x)}{(2 \pi)^{{d_v}/2}}\exp \left( -\frac{\left |v\right |^2}{2} \right), \quad \forall (x,v) \in \TT^{d_x} \times \RR^{d_v}, \\
    & \omega(0, x) = \sqrt{\frac{2}{{d_v}\, \rho_0(x)}},
  \end{aligned} \right.
\end{equation}
with ${d_x}={d_v}=1$, $\rho_{0}(x)= 1+0.1 \cos(\pi x)$ and $\TT=[0,L]$. The fact that we used an unidimensional velocity space is relevant, because we consider inelastic collisions.  Indeed, for this type of collision, the  dissipation of kinetic energy is enough to insure that the microscopic dynamics is nontrivial (which is of course not the case with elastic collisions).

\vspace{0.2cm}

\paragraph{\bf Description  of the results.} 

We start with a simple one dimensional test intended to show that the uncoupled model is a good approximation of equation \eqref{eqCollision} in the space inhomogeneous case. 

We take $L = 1$ and $N_x=25$ points in space variable, a rescaled velocity box $\mathcal{D}_V$ of size $V=10$ and $N_v=32$ half Fourier modes in rescaled variables. We also perform computations in classical velocity variables, with a velocity box of size $V=8$ and $N_v = 128$ half modes, in order to achieve a good accuracy and avoid oscillations in short time. The initial condition is given by \eqref{eqCItrendToEqu}.

\begin{figure}
\begin{center}
\includegraphics{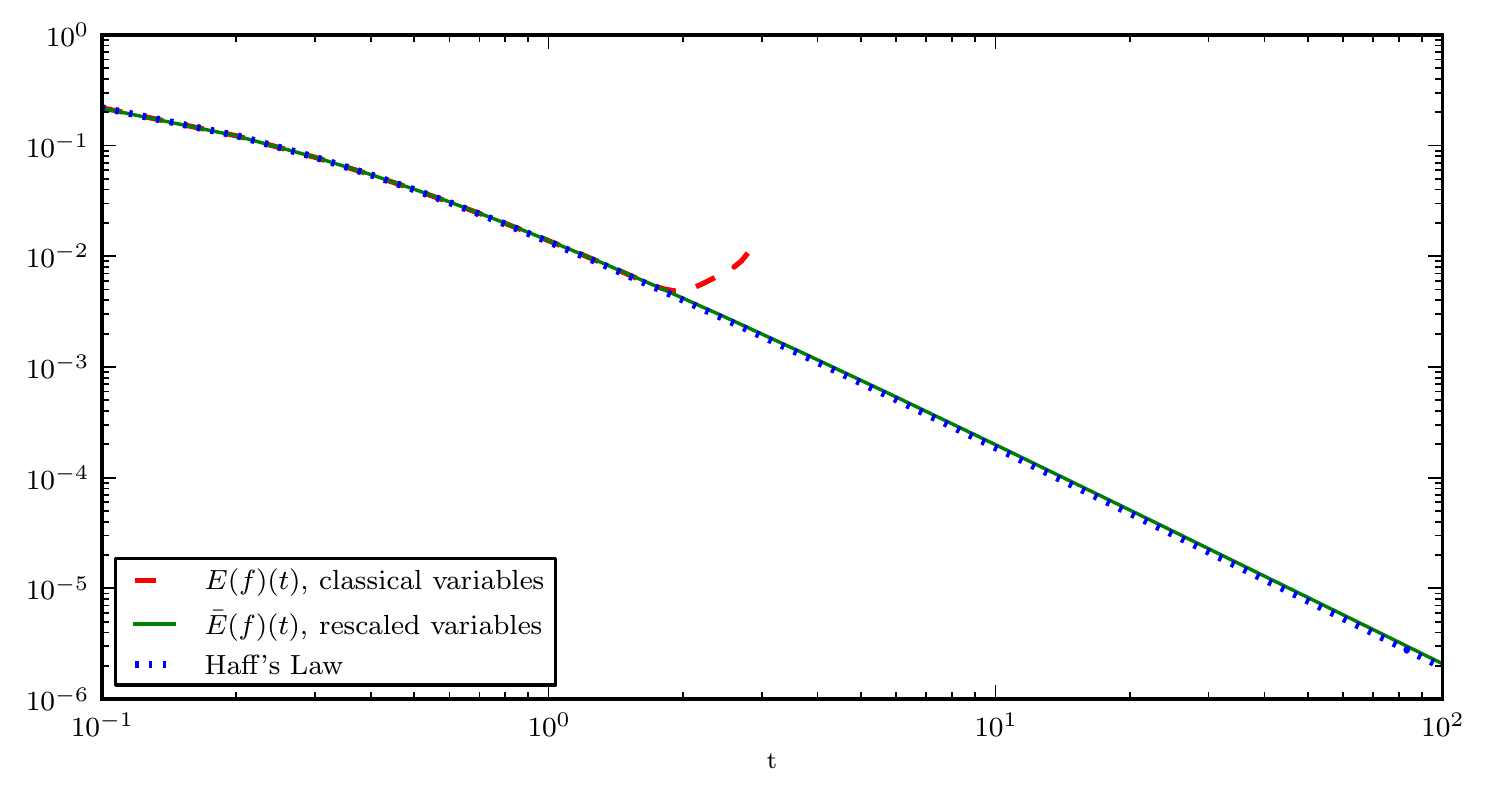}
\caption{Test 2 -- Comparison of the global inhomogeneous cooling process in original and scaled variables, with $e=0.9$.}
\label{figCompHaff}
\end{center}
\end{figure}

We compare the evolution of the global kinetic energy ${E}(t)$ of both simulations. We also plot the theoretical asymptotic behavior of the kinetic energy (namely, the Haff's Law \cite{haff:1983}), given formally by
\begin{equation*}
E(t) \simeq \frac{E(0)}{(1 \,+\, \sqrt{2} K_1 \,t)^2}.
\end{equation*}
We can see in Figure \ref{figCompHaff} that the Haff's Law is very well captured by the rescaling method. This is also the case for short time in the classical method, but due to the fast apparition of spurious oscillations (even with $N_v=128$ half modes), a strong divergence from the expected global behavior occurs in this case. Once more, the efficiency of the rescaling method for the long time behavior of the solution becomes clear, since the Haff's Law is perfectly described.

%\vspace{0.2cm}
%\noindent
%\emph{Test 2. Oscillations in the trend to equilibrium.}

\begin{figure}
  \begin{center}
    \begin{tabular}{cc}
      \includegraphics[scale=1]{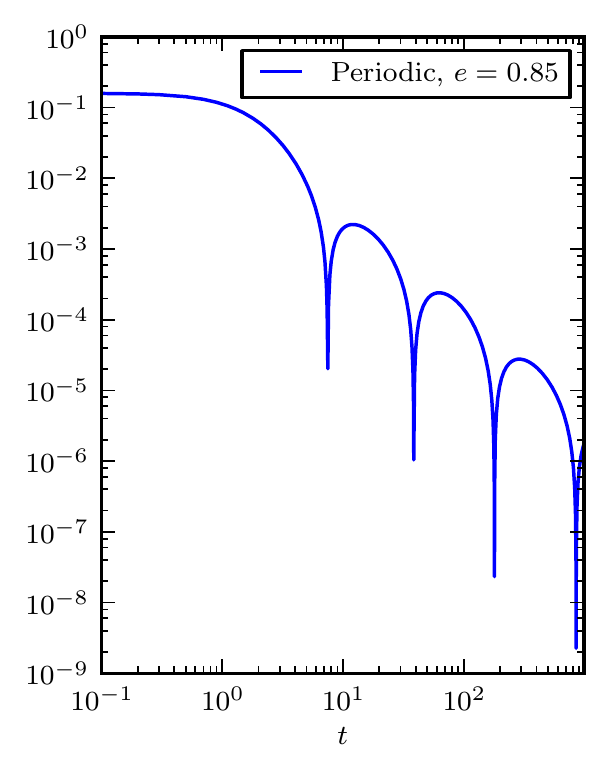} & \includegraphics[scale=1]{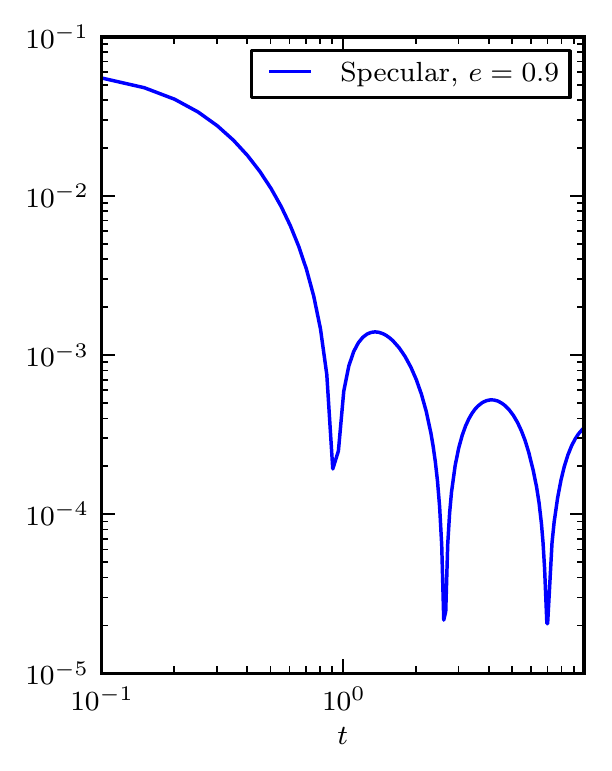}
    \end{tabular}
    \caption{Test 2 -- Oscillations of the quantity $\mathcal L(t)$ for periodic and specular boundary conditions. }
    \label{figretourEquBC}
  \end{center}
\end{figure}

\begin{figure}
	\begin{center}
		\includegraphics{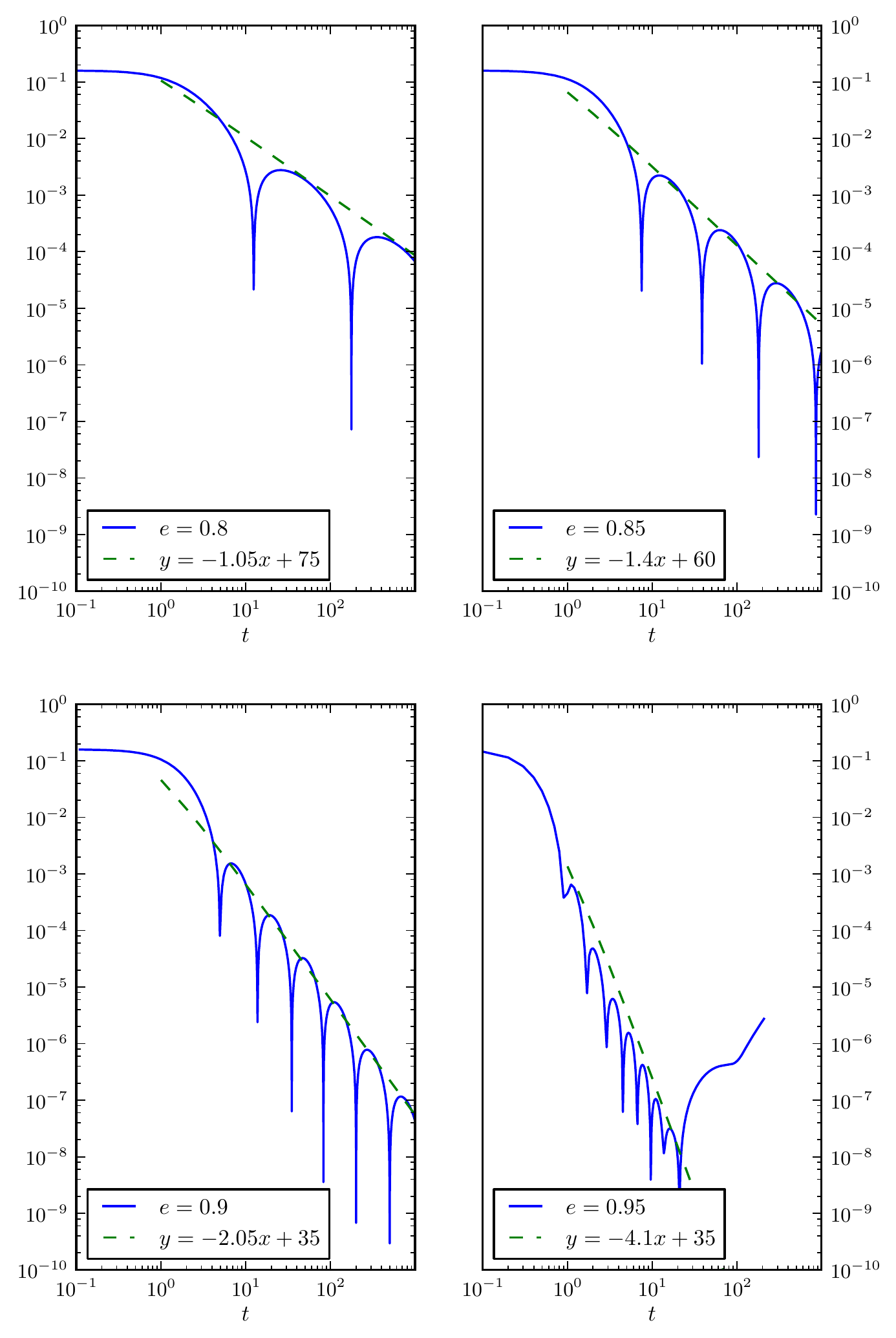}
		\caption{Test 2 -- Time evolution of  $\mathcal L(t)$ (deviation from equilibrium of the local mass), for different values of the restitution coefficient $e$.}
		\label{figretourEquMulti}
	\end{center}
\end{figure}

Let us investigate more precisely the long time behavior of the numerical solution obtained from the rescaled method. In Figures \ref{figretourEquBC} and \ref{figretourEquMulti}, we  present the time evolution in $\log$-$\log$ scale of 
\[\mathcal{L}(t) := \|\rho(t, \cdot) - 1 \|_{L^1}\]
for an initial condition given by \eqref{eqCItrendToEqu}.
This quantity describes the relaxation towards global equilibrium of the local mass $\rho(t,x)$ in $L^1$ norm.
 
Numerical investigations on the elastic Boltzmann equation showed some evidences of oscillations during the relaxation process, for example in the local relative entropy \cite{filbet:2007}, as it was conjectured by L. Desvillettes and C. Villani in \cite{desvillettes:2005}. In the inelastic case, the situation is different. On the one hand the solution  $f(t)$ does not converge to a smooth steady state since $f_\infty$ is the Dirac mass \eqref{maxwglob}, but hopefully our rescaling method allows to follow the concentration in velocity. On the other hand, on the evolution of $\mathcal{L}(t)$ we  still observe the apparition of oscillations, but now the oscillation frequency is constant with respect to the logarithmic time scale, which is a big difference with the elastic case (see Figure~\ref{figretourEquBC}).

  We also observe in Figure \ref{figretourEquBC} that the frequency of oscillation is different whether we consider specular or periodic boundary conditions. This is due to the fact that this frequency depends on the size of the box, and that specular conditions on a box are equivalent to periodic conditions on $2^{d_v}$ copies of this box, by symmetrization.

Finally, we point out that the scheme breaks down when $\mathcal L(t)$ becomes too small (this is particularly clear for the case $ \ve = 0.05$ in Figure \ref{figretourEquMulti}). This can be avoided by refining the space and velocity grids.
				
\vspace{0.2cm}

\paragraph{\bf Interpretation of the results.} Here, we performed simulations on the full inelastic Boltzmann equation in a
simplified geometry (one dimension of space and velocity, periodic and specular boundary conditions, fixed Knudsen number) with the spectral method to observe the evolution of the quantity $\mathcal{L}(t)$ and to check numerically if such oscillations occur.  It is in fact possible to give a simple interpretation of  these oscillations thanks to the work \cite{ellis:1975}.  Since this effect is observed near the global equilibrium, which is homogeneous in space, one can consider that the space inhomogeneous granular gases equation is a small perturbation of the homogeneous equation, which exhibits a self-similar behavior and the Haff's Law. Therefore,  we rescale this equation as
\begin{equation}
\label{defScalHomog}
\left\{ \begin{array}{l} 
f(t, x, v) = V(t)^{d_v} h(s(t),x, V(t) \,v),
 \\
\,
\\
V(t)= 1 \,+\, (1-e) \, t, 
\\
\,
\\
\displaystyle s(t)= \frac{\log(1 + (1-e) \, t)}{(1-e)}.
\end{array}\right.
\end{equation}
Then, one can show after straightforward computations that the distribution $h$ is solution to the granular Boltzmann equation with an anti-drift term in the usual set (see the paper \cite{Mischler:20062} and the references therein) of self similar variables $(s, x, \nu) := (s(t), x, V(t)\,v)$
\begin{equation} 
\label{eqdriftInh}
\frac{\partial h}{\partial s} \,+\, \nu \cdot \nabla_x h  \,+\, (1-e) \, \nabla_\nu (\nu h)   \,=\, \frac{\Q_{\mathcal I}(h, h)}{\ve}.
\end{equation}

Due to the convergence toward an equilibrium distribution $H_e$, we can replace the collision operator by its linearization $\mathcal L_e$ around $ H_e$ in \eqref{eqdriftInh}. The initial condition \eqref{eqCItrendToEqu} actually amounts to study thanks to $\mathcal L(t)$ the long time behavior of the first spatial Fourier mode of $h$. This behavior is formally given by the one of $t \mapsto e^{\lambda_j t}$ when $t \to \infty$ for any $j\in\{0,\ldots,{d_v}+1\}$, where $\lambda_j$ is one of the ${d_v}+2$ hydrodynamic eigenvalues of the linearized operator $\mathcal L_e$. 
				 
At least for weak inelasticity, we have seen in Subsection~\ref{subGranular} that the inelastic collision operator can be understood as a perturbation of the elastic Boltzmann operator. In \cite{Rey:2012EllisPinsky}, we used this expansion  and results from \cite{mischler:20091} to perform a spectral analysis of $\mathcal L_e$, as it has already been done for the elastic Boltzmann operator in \cite{ellis:1975}. This allows to compute the first terms of the Taylor expansion of the hydrodynamic eigenvalues $\lambda_j$.
We give this expansion with respect to $\eta = 2 \pi/L$ and $1-e$:
\begin{equation}
\label{eigExpans}
\lambda_j \,=\, i \,\eta\, \lambda_j^{(1)} \,+\, \eta^2 \,\lambda_j^{(2)} \,-\, (1-e)  \,+\, \mathcal O(\eta^3) \,+\, \mathcal O((1-e)^2),
\end{equation}
where $\lambda_j^{(1)}$ is the first order expansion in $\eta$ of the elastic hydrodynamic eigenvalues:
\begin{equation*}
\lambda_j^{(1)} \,=\, \left\{
\begin{array}{ll}
\displaystyle \pm \sqrt{1+\frac{2}{{d_v}}}, & j \in \{0,1 \}, 
\\
0,& j\in\{2,...,{d_v}+1\}. 
\end{array} \right .
\end{equation*}	
The second order expansion in $\eta$ is given by
\begin{equation*}
\lambda_j^{(2)} = \left\{\begin{aligned}
& -\frac{\kappa}{{d_v}+2} - \frac{\mu}{2}, && j \in \{0,1 \}, \\
& -\left( 1 + \frac{2}{{d_v}}\right )\kappa, && j=2, \\
& -\frac{{d_v}}{2({d_v}-1)} \, \mu, &&  j\in\{3,...,{d_v}+1\},
\end{aligned} \right .
\end{equation*}
with $\mu$ the viscosity coefficient and $\kappa$ the heat conductivity of the Navier-Stokes limit of the elastic Boltzmann equation.

The expansion \eqref{eigExpans} can explain the oscillations and the damping observed in Figure~\ref{figretourEquMulti}. 
Indeed, the first order term in $\eta$ is imaginary and will generate oscillations during the return toward equilibrium, of period $\eta \sqrt{1+2/{d_v}}$. 
Then, the period of oscillation $\zeta$ will be proportional to $L^{-1}$, as in \cite{filbet:2007}. 
In our case, the time has not been scaled, and then according to the definition \eqref{defScalHomog} of $s(t)$, the oscillations will appears in a logarithmic time.
The second order term in $\eta$ and the term in $(1-e)$ are nonpositive and then will exhibit a damping effect of magnitude $\eta^2 \alpha + (1-e) \beta$ for $\alpha$, $\beta$ positive, according to \eqref{eigExpans}. 
We can see in Table \ref{tabNum} that this effect can be observed numerically: when the restitution coefficient vary, the quantity $ (1-e) \beta$ remains constant, when measured in the logarithmic scale.

\begin{table}
	\begin{center} \begin{tabular}{|c|c|c|}
	\hline Restitution coefficient $e$ & Damping rate $\beta$ & $ - (1-e) \beta$ \\
	\hline 0.8 & -1.05 & 0.21 \\
	\hline 0.85 & -1.4 & 0.21 \\
	\hline 0.9 & -2.05 & 0.205 \\
	\hline 0.95 & -4.1 & 0.205 \\
	\hline
	\end{tabular} \end{center}
	\caption{Test 2 -- Influence of the restitution coefficient $e$ on the damping rate in the time evolution of the quantity $\mathcal L$.} 
	\label{tabNum}
\end{table}	

\subsection{Test 3: Inelastic Shear Flow Problem}

	The uniform shear flow problem is perhaps one of the most widely studied nonequilibrum problem for granular gases, both experimentally and theoretically \cite{Astillero:2005}. It can be described as follows: consider a gas enclosed between two infinite planes located at $y = -1/2$ and $y=1/2$, in opposite relative motion with respective velocities $-1/2$ and $1/2$ in the $x$-direction. The planes are not described as realistic thermalized, reflecting walls (\ie Maxwell boundary conditions) but as Lee-Edwards boundaries: when a particle of velocity $v$ hits one of the planes, it is reemitted through the opposite one with a velocity $v'$ such that the relative velocity of the particle with respect to the plane is preserved before and after this exchange. In consequence, the energy of the system increases during this process.
	
	Thus, the particles gain energy at the boundary while other particles lose energy due to the inelasticity of the collisions. This problem has been solved numerically using Direct Simulation Monte Carlo methods in \cite{Astillero:2005,Ruiz-Montero:2010}. We shall consider it using the full deterministic inelastic Boltzmann equation. 
	The gas is taken initially at thermal equilibrium in the two dimensional infinite tube
	\begin{equation*}
	  f_{0}(x,v) = \frac{1}{2 \pi}\exp \left( -\frac{\left |v\right |^2}{2} \right), \quad \forall \, x \in [-1/2, 1/2] \times \RR, \quad \forall \, v \in \RR^2.
	\end{equation*}
	The evolution of the system will be induced by the inelastic collisions as well as the flow created by the energy input at the boundaries.

  \vspace{0.2cm}
  
	\paragraph{\bf Description  of the results.} 
		
  The simulations are made thanks to the rescaled model \eqref{eqGgran}, with the granular gases collision operator $\Q_\mathcal{I}$, on a simplified geometry: one dimension of space and velocity with $\Omega = [-1/2, 1/2]$. 
  The shear flow problem is mimicked by using diffusive boundary conditions: we impose  constant temperature $T_w = 1$ and velocity $u_w = 0$ on the boundary $x = \pm \, 1/2$. The restitution coefficient of the inelastic collisions is $e=0.9$. We take $N_x = 100$ and $150$ points for the space grid and $N_v = 32$ half Fourier modes  on a rescaled velocity box $\mathcal D_V = [-V,V]$ of size $V=12$.
  We compare these results with computations made on the model \eqref{eqCollision} in classical variables. We take in this case $N_x = 100$ spatial points and $N_v = 32$ half Fourier modes  on a velocity box $\mathcal D_V$ with $V=8$.

	\begin{figure}
		\begin{tabular}{c}
		  Normalized density \\
		  \includegraphics{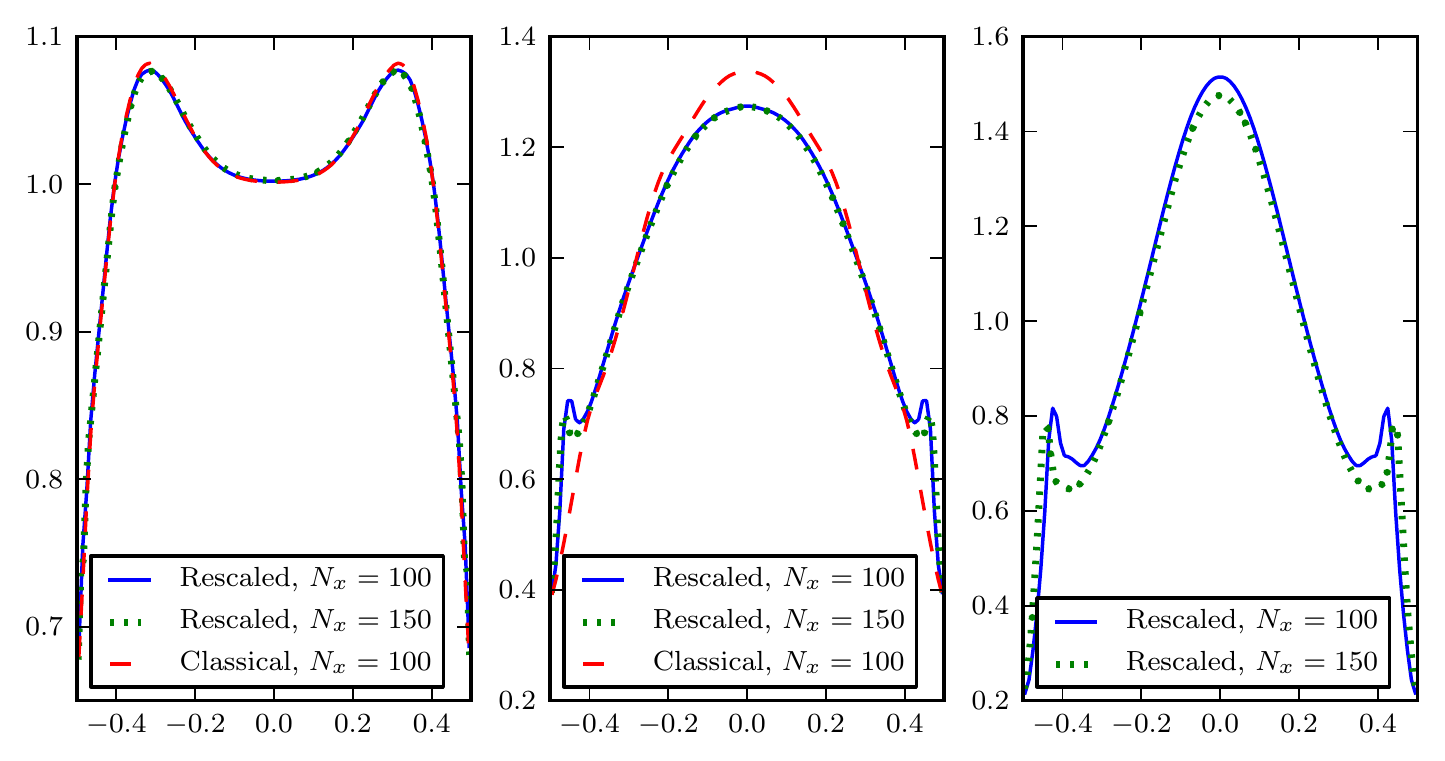} \\
		  Normalized temperature \\ 
		  \includegraphics{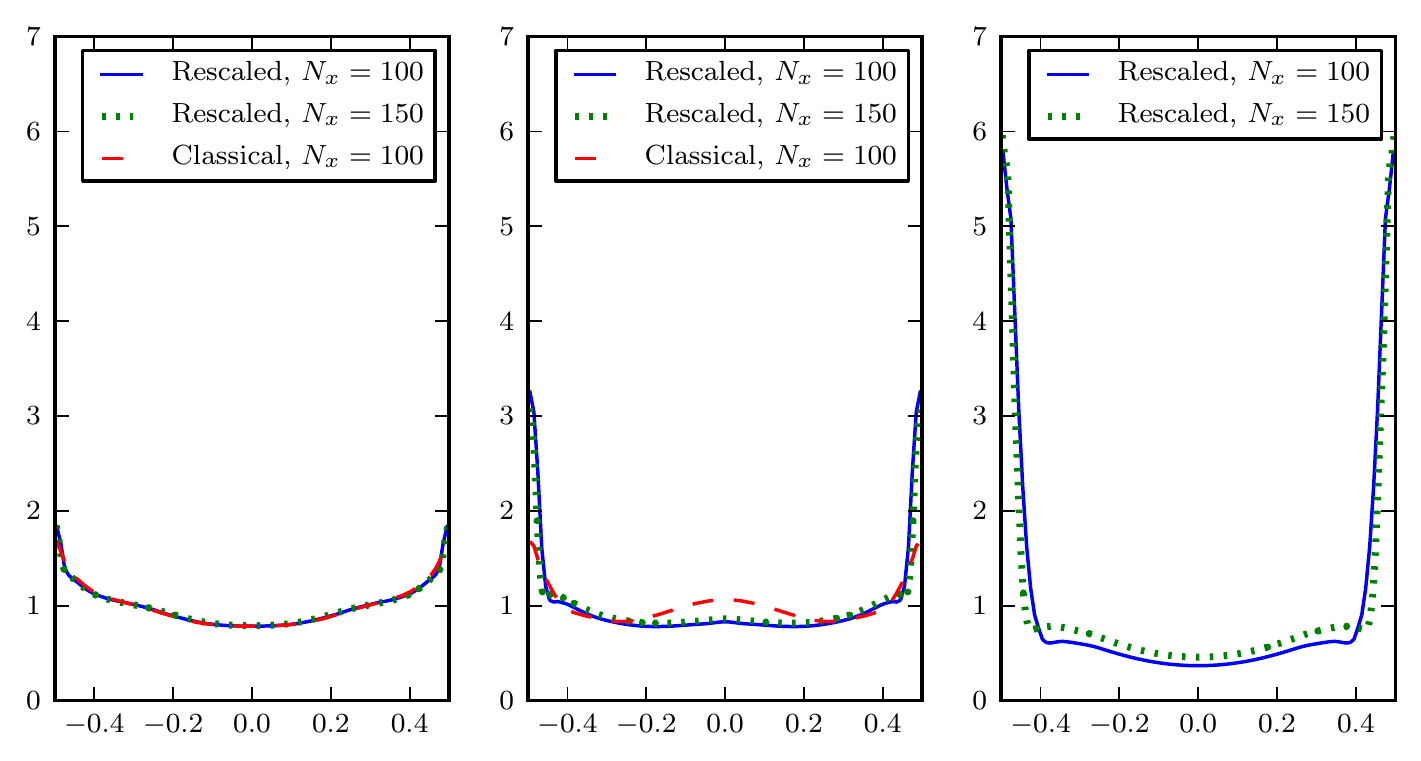}	
		\end{tabular}
		\caption{Test 3 -- Normalized density $\rho(t,x)/ \overline \rho$ and temperature $T(t,x)/\overline T(t)$ of the shear flow problem computed in classical and rescaled variables, at times $t=0.25$, $0.75$ and $2$.}
		\label{figShearFlow}
	\end{figure}
	 
	  We present in Figure \ref{figShearFlow} the evolution of the normalized density $\rho(t,x)/ \overline \rho$ and temperature $T(t,x)/\overline T(t)$, where $\overline \rho$ and $\overline T$ are the space averaged quantities
	  \[ \overline{T}(t) := \frac{1}{|\Omega|} \int_{\Omega} T(t,x) \, dx, \quad \overline{\rho} := \frac{1}{|\Omega|} \int_{\Omega} \rho(t,x) \, dx. \] 
    On the one hand, we observe that the computation in classical and rescaled variables agree very well for short time. Then, the oscillations due to the concentration in velocity variable for the classical case quickly break the result of this computation. On the other hand, the solution obtained in rescaled variables can be computed for large time, thanks to the change of scales induced by the choice \eqref{eqWclose} of $\omega$. We also notice the stability of the scheme with respect of the space grid.
  
  The short time behavior of this quantity is close to the one presented in \cite{Astillero:2005} using DSMC simulations on the more realistic $2d_x \times 2d_v$ geometry. Indeed, we can see that the thermalized walls heat the system at the boundary and the collisions cool it inside the domain. This induces very strong gradients near the boundaries for the macroscopic quantities. When time evolves, the temperature will eventually become uniform inside the domain with bigger and bigger gradients near the boundary. 
  This large time behavior is different of the one described in \cite{Astillero:2005} where the temperature reaches a space uniform steady state. This might be due to the use of a simplified geometry. Indeed, in their case, the heating of the boundary induced by the Lee-Edwards conditions will stop because of the uniformization of the gas.

\subsection{Test 4: Cluster Formation in Granular Gases}
\label{subFPBGK}

	We consider now the kinetic equation $\eqref{eqCollision}$ with the simplified granular gases operator $\Q_\mathcal{S}$ in two-dimension in space on the torus $\TT^2 = [0,1]^2$, and  three-dimension in velocity
	\[
	\frac{\partial f}{\partial t} + v \cdot \nabla_x f = \frac{\Q_{\mathcal{S}}(f)}{\ve}, \quad \forall \, (x,v) \in \TT^2 \times \RR^3
	\]
	with periodic boundary conditions in $x$. We shall present some numerical evidences of \emph{clustering} for this simplified granular gases model.
	
	The problem of cluster formation in a force-free granular media is a particularly difficult test case of numerical methods for granular gases simulations. It has been shown to occur numerically in both microscopic and macroscopic simulations.
	More precisely, in the classical Molecular Dynamics (MD) simulations of~\cite{Goldhirsch:1993}, clusters arise from an initially homogeneous gas of particles. It was also observed using MD simulations in \cite{Poschel:2005} that clustering is a transient phenomenon: an initially homogeneous system of particles interacting by inelastic collisions will tend to form clusters and then will eventually come back after some time to the homogeneous cooling state.
	Clusters were also observed using deterministic  numerical simulations for the macroscopic quantities, thanks to the energy dissipative Navier-Stokes model of~\cite{carrillo:2008granular}.
	
	We will show that this phenomenon also occurs at the mesoscopic scale, namely using a kinetic description.
  The initial condition for the problem is taken as a Maxwellian distribution of homogeneous, isotropic momentum and temperature and uniform density:
  \begin{equation*}
    f_{0}(x,v) = \frac{U(x)}{(2 \pi)^{{d_v}/2}}\exp \left( -\frac{\left |v\right |^2}{2} \right), \quad \forall \, (x,v) \in \TT^{2} \times \RR^{3}
  \end{equation*}
  with $U$ an uniformly distributed random variable over $[0,1]^2$.
	
  \vspace{0.2cm}
  
	\paragraph{\bf Description  of the results.} 
	
	The simulations are made using the simplified granular gases operator $\Q_\mathcal{S}$ presented in Subsection~\ref{subFPBGK} in the rescaled model \eqref{eqGcons}-\eqref{eqSysConsUncoupled}. The macroscopic scaling function $\omega$, given by \eqref{eqWclose}, is computed thanks to the energy dissipative Euler system \eqref{eqSysCons}. Due to the correct energy dissipation, this model mimics the behavior of the kinetic granular gases equation and is then consistent with the MD simulations of~\cite{Goldhirsch:1993}. Moreover, the ES-BGK operator used in the definition of $\Q_\mathcal{S}$ gives the good Prandtl number for the Navier-Stokes limit of the model, which is also consistent with the macroscopic model of~\cite{carrillo:2008granular}.
	
	\vspace{0.2cm}
	\noindent
		
\begin{figure}
	  \begin{tabular}{c}
	    Kinetic density \\
	    \includegraphics{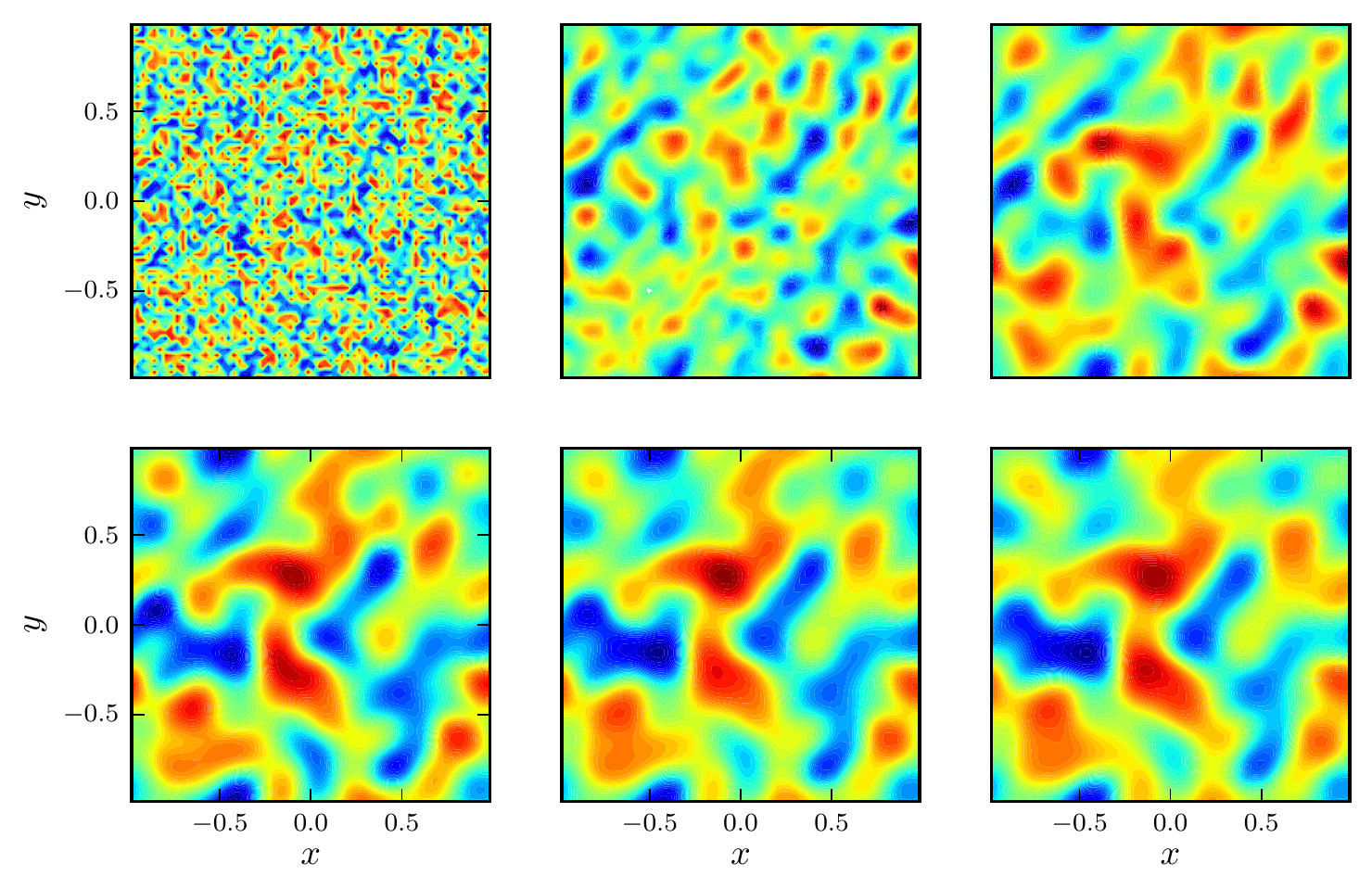} \\
	    Hydrodynamic density \\
	    \includegraphics{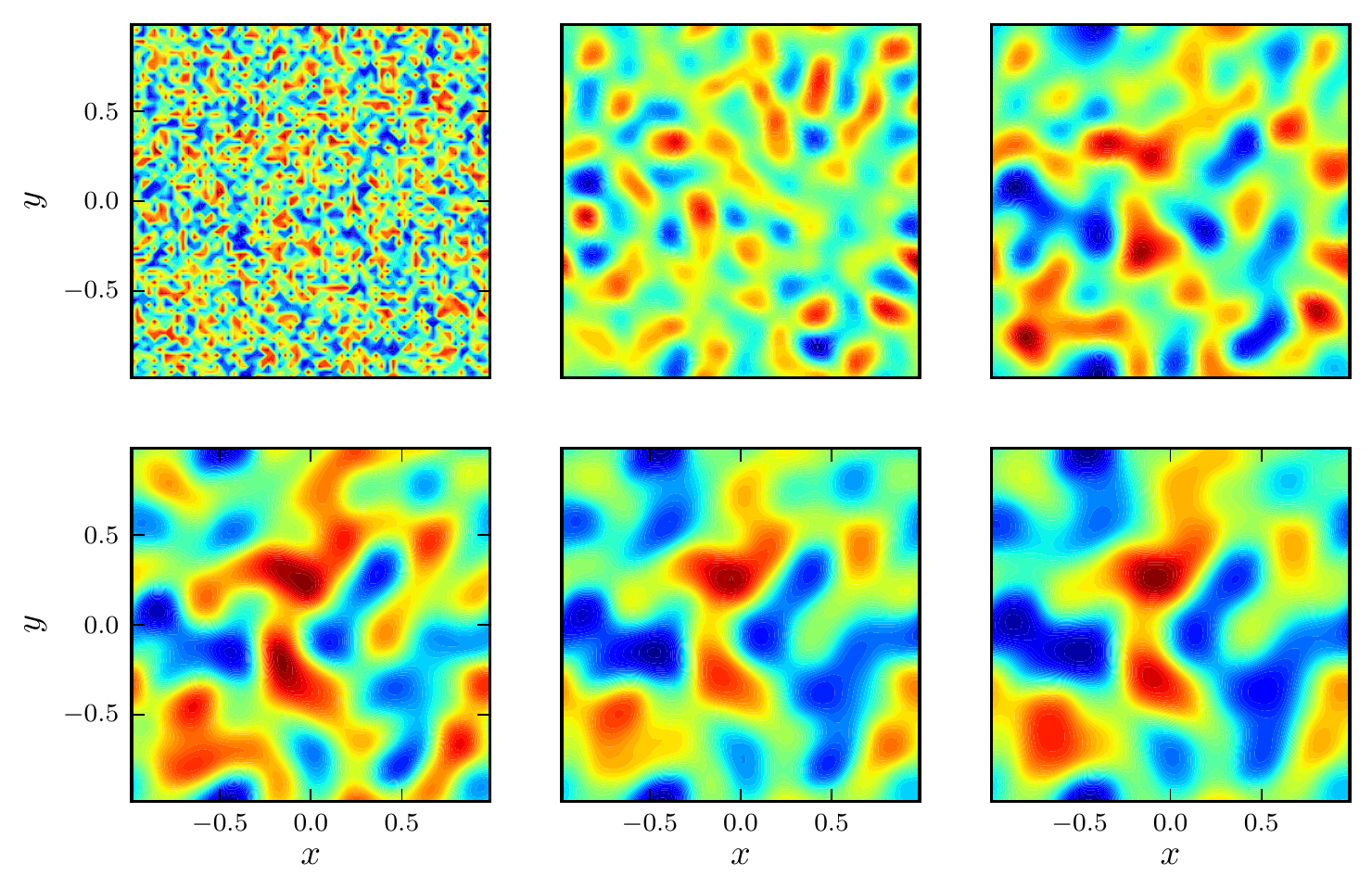}
	  \end{tabular}
	   \caption{Test 4 with $e=0.99$ -- Contour plot of the density given by the simplified granular gases model (top) and the dissipative Euler system (bottom), at times $t=0$, $0.25$, $1$, $5$, $10$ and $30$. High values are displayed in red, low in blue.}
	   \label{figClusterQuasiElastic}
	\end{figure}

	We take $N_x = 64$ points in each direction of space on the torus and a velocity box of size $V=7$ with $24$ points in the directions $v_x$ and $v_y$ and $8$ in the direction $v_z$. The restitution coefficient is $e = 0.99$ (and then $\ve = 0.01$ thanks to the weak inelasticity asumption \eqref{quasi}).
	We show in Figure~\ref{figClusterQuasiElastic} the density contour plot given by the kinetic model and by the dissipative Euler system~\eqref{eqSysCons} used to compute the scaling function. We see that both models give close results for short time, which is consistent with the fact that the hydrodynamic limit of the kinetic equation \eqref{eqCollision} with the Fokker-Planck operator \eqref{eqOpFPBGK} is given at the first order in $\ve$ by Euler system~\eqref{eqSysConsUncoupled}. Then, for larger time, clusters start to appears even with a coarse grid in space.
	It is interesting to notice that the clusters remains bounded through time, even without taking into account the dense gas effects through a ``pair correlation function'' depending of a maximal packing fraction, as in \cite{carrillo:2008granular}.
	 
\begin{figure}
	  \begin{tabular}{c}
	    Kinetic density \\
	    \includegraphics{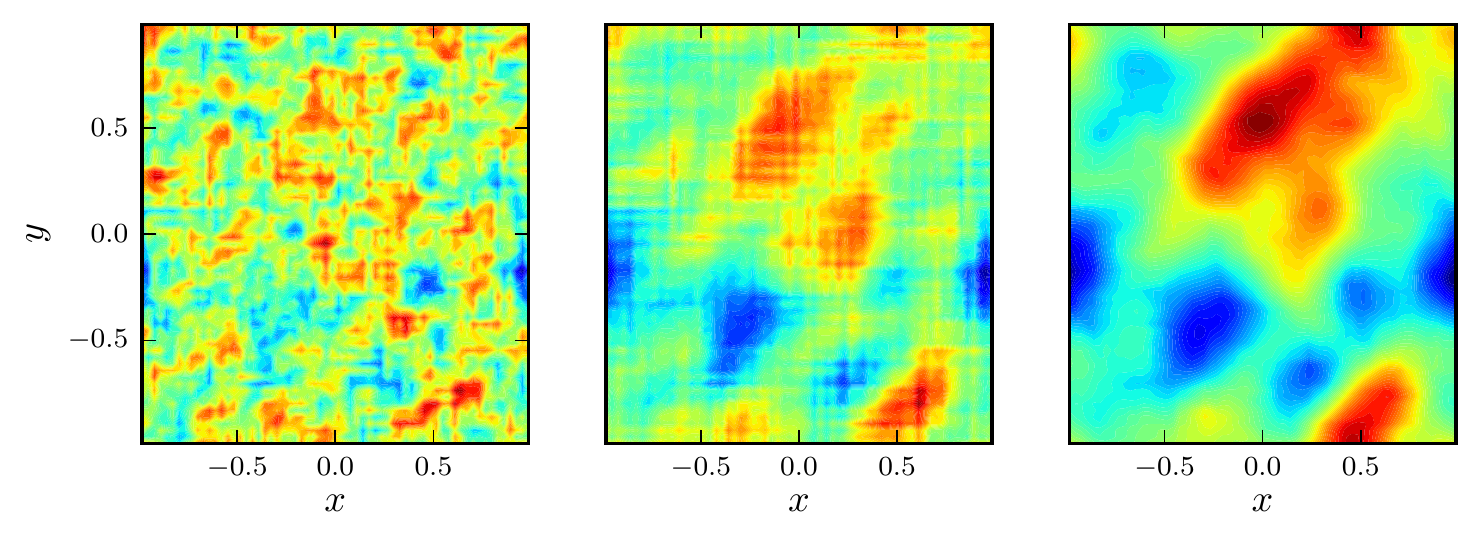} \\
	    Hydrodynamic density \\
	    \includegraphics{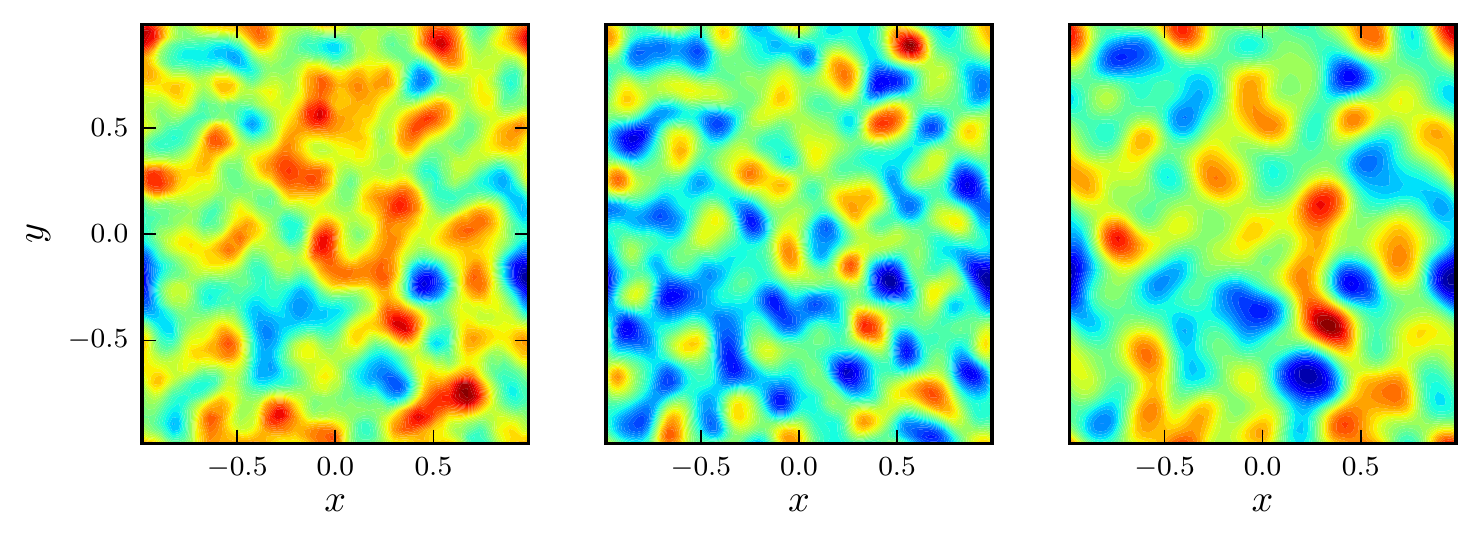}
	  \end{tabular}
	   \caption{Test 4 with $e=0.75$  -- Contour plot of the density given by the simplified granular gases model (top) and the dissipative Euler system (bottom), at time $t=0.05$,  $0.5$ and $1$. }
	   \label{figClusterNotQuasiElastic}
	\end{figure}

\begin{figure}
	  \begin{tabular}{cc}
	    Kinetic temperature                             &  Hydrodynamic temperature \\
	    \includegraphics{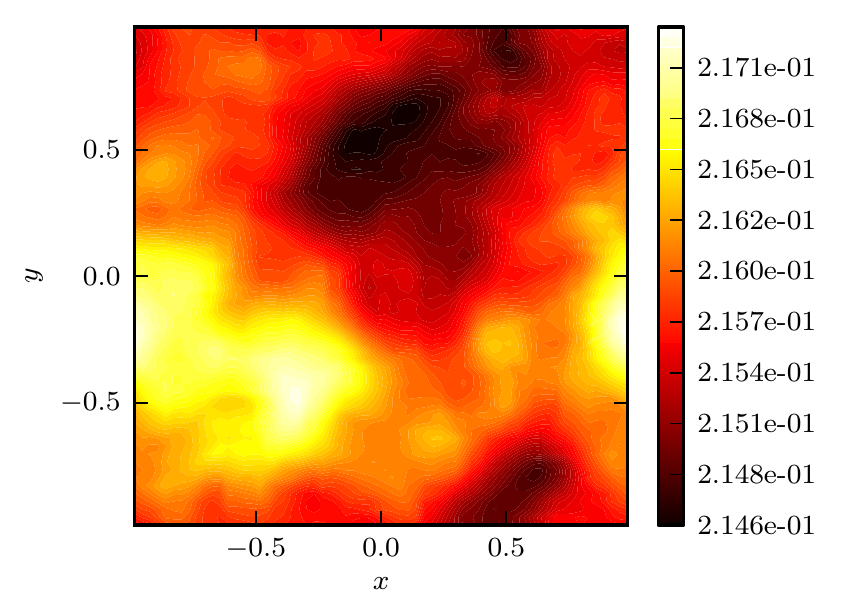} & \includegraphics{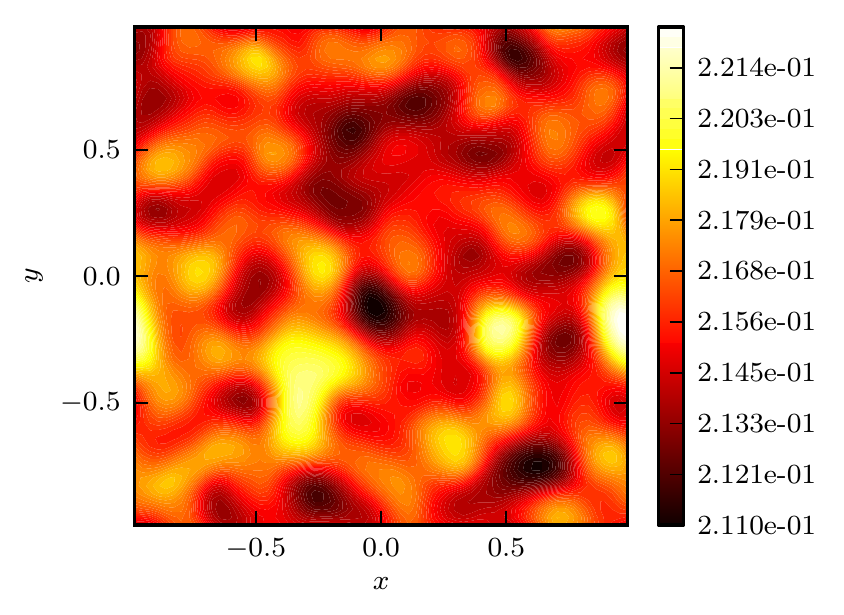}
	  \end{tabular}
	   \caption{Test 4 with $e=0.75$  -- Contour plot of the temperature given by the simplified granular gases model (top) and the dissipative Euler system (bottom), at time  $t=1$. }
	   \label{figClusterNotQuasiElasticTemp}
	\end{figure}

  We then present a test where the restitution coefficient is taken as $e = 0.75$, with the same  computational grid than before. This example corresponds to a case where the hydrodynamic closure we used \emph{does not} apply to the model, because $\ve$ is ``large''.
  Hence, we observe in Figure \ref{figClusterNotQuasiElastic} that the hydrodynamic density does not correspond at all with the kinetic one (even in short time). Nevertheless, our method works well even in this case, because the temperature given by the macroscopic model is of the good order of magnitude compared to the one we obtain with the kinetic model, as shown in Figure \ref{figClusterNotQuasiElasticTemp} (there is an agreement on the temperature up to $3$ digits). Hence, the scale of the distribution function in the velocity space is well resolved by the kinetic scheme, even if the temperature is very low (the order of magnitude at time $t=1$ is about $0.2$) and the velocity grid very coarse.
  It then confirms that more than the accurate resolution of the hydrodynamic fields, it is the scale of their temperature which matters in our method.

\subsection{Test 5: Riemann Problem for the Elastic Boltzmann Equation}

  To finish, we consider  the full kinetic equation \eqref{eqCollision} with the ES-BGK operator for elastic interactions (as presented in Section \ref{subBGK} or in \cite{AndriesPerthame:2000}), in two dimensions of space and three dimensions of velocity, for an infinite domain in space:
  \[
	  \frac{\partial f}{\partial t} + v \cdot \nabla_x f = \frac{ \tau(\rho_f) \,\left( \G[f] \,-\, f \right)}{\ve}, \quad \forall \, (x,v) \in \RR^2 \times \RR^3,
	\]
	with $\tau(\rho) = 0.925 \, \rho \, \pi /2$.
	Even if the ES-BGK operator (as well as the elastic Boltzmann operator) does not exhibit concentration nor spreading in the velocity variables, our rescaling velocity method can be of great interest if one wants to simulate some problems with a huge difference of magnitude in temperature.
	
\begin{figure}
    \begin{center}
      \includegraphics[scale=1]{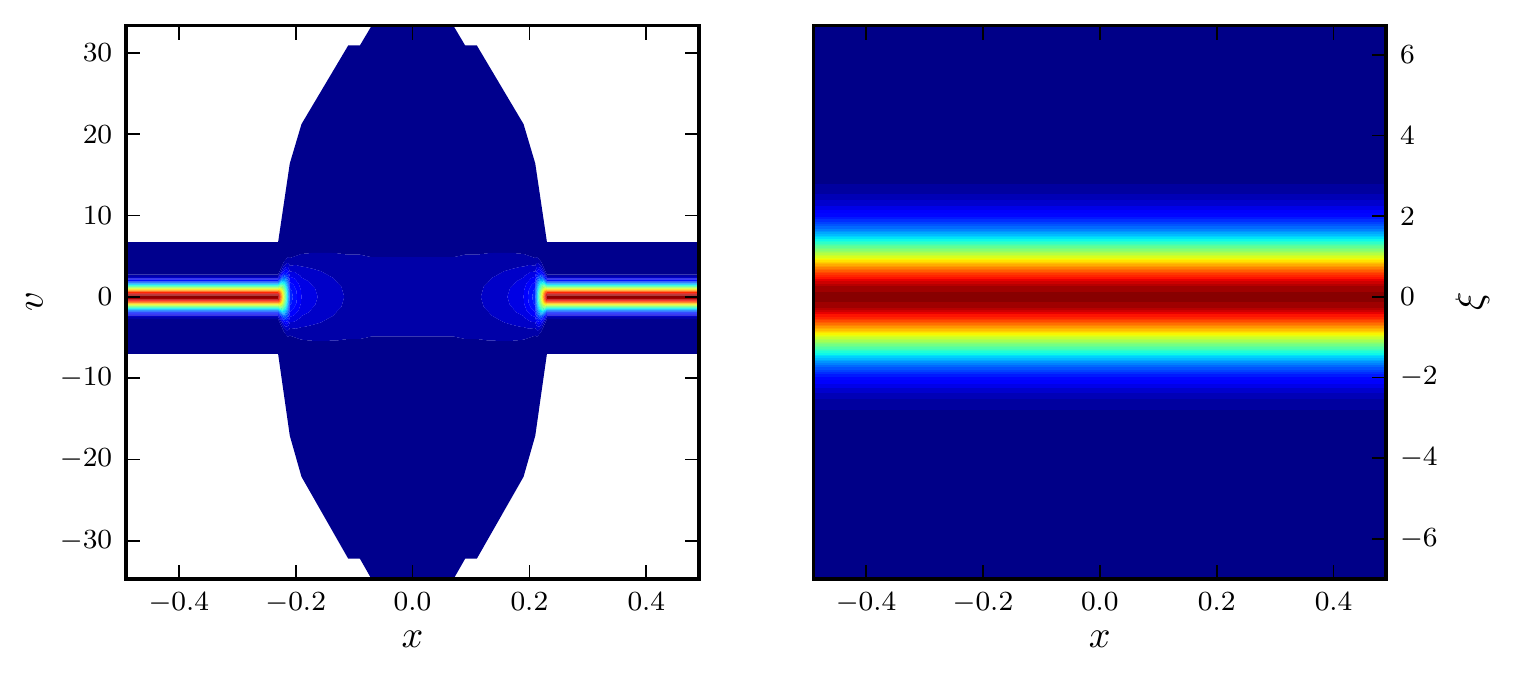}
      \caption{Test 5 -- First marginal in space and velocity of the initial condition for a shock in temperature, in both classical and scaled velocity variables}
      \label{figCompShockFG}
    \end{center}
  \end{figure}

	Let us then consider a gas at rest in an infinite space, which is suddenly heated at one particular point. The change in temperature  will induce a shock, which will propagates in the gas. If this change is big enough, the particle distribution function will be very concentrated in velocity where the gas is cold, and very spread where it is hot.
	We will then take as an initial condition the following centered Gaussian distribution
	\begin{equation*}
    f_{0}(x,v) = \frac{1}{(2 \, \pi \, T_0(x))^{d/2}}\exp \left( -\frac{\left |v\right |^2}{2 \,T_0(x)} \right), \quad \forall \, (x,v) \in \RR^2 \times \RR^3,
  \end{equation*}
  where the initial temperature $T_0$ verify
  \begin{equation*}
    T_0(x) = \left \{ 
      \begin{aligned} 
        1   & \text{, if } |x| \leq 0.5, \\ 
        100 & \text{, if } |x| > 0.5. 
      \end{aligned} \right. 
  \end{equation*}
  The rescaled initial condition is then set as $g_0(\cdot, \xi) = \omega(0, \cdot)^d f_0\left (\cdot, \omega(0, \cdot) \, \xi\right )$ for all $\xi \in \RR^3$, with
	\begin{equation*}
	  \omega(0, x) = \sqrt{T_0(x)}.
	\end{equation*}
  We presented in Figure \ref{figCompShockFG} the first marginal of this initial distribution, in both scaled and classical velocity variables. We can see that if one wants to solve the kinetic equation \eqref{eqCollision} in classical variables for such a distribution, one has to take a very large velocity space to capture the hot part of the gas, but the grid also has to be very thin in order to resolve the cold part, which is much more concentrated near the origin.
  Then, the number of points in the velocity space (which is two dimensional) has to be huge.
  
  With the rescaling velocity method, this is not a problem anymore, because as one can see in the right hand side part of Figure \ref{figCompShockFG}, the initial condition in the rescaled space is a simple reduced Maxwellian distribution, uniform in velocity. All the information on the initial distribution will be contained in the scaling function $\omega$. This function will be evolved thanks to the Euler system of perfect gases \eqref{eq:macro:elastic}.
  
  \vspace{0.2cm}
  
	\paragraph{\bf Description  of the results.} 
	
	  We take $N_x = 64$ points in each direction of space (with absorbing conditions to mimic an infinite domain) and and a velocity box of size $V=7$ with $24$ points in the directions $v_x$ and $v_y$ and $8$ in the direction $v_z$. We compute our numerical solution in an hydrodynamic regime, in order to compare it with the compressible Euler system \eqref{eq:macro:elastic}. Hence, the Knudsen number is taken as $\ve = 0.05$ in order to have close enough macroscopic solutions.
		
\begin{figure}
	  \begin{center}
	    Kinetic model
		  \begin{tabular}{cc}
		    Density                                         & Velocity ($x-$component)                       \\
		    \includegraphics{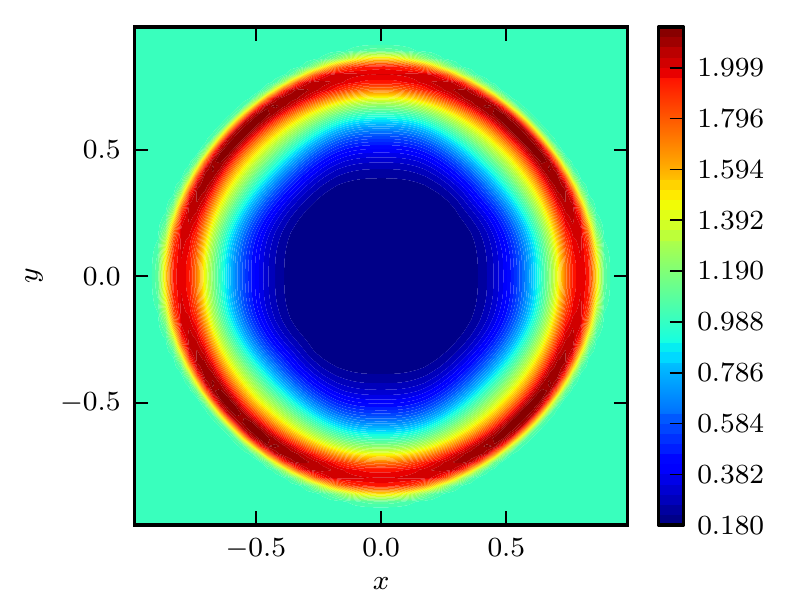} & \includegraphics{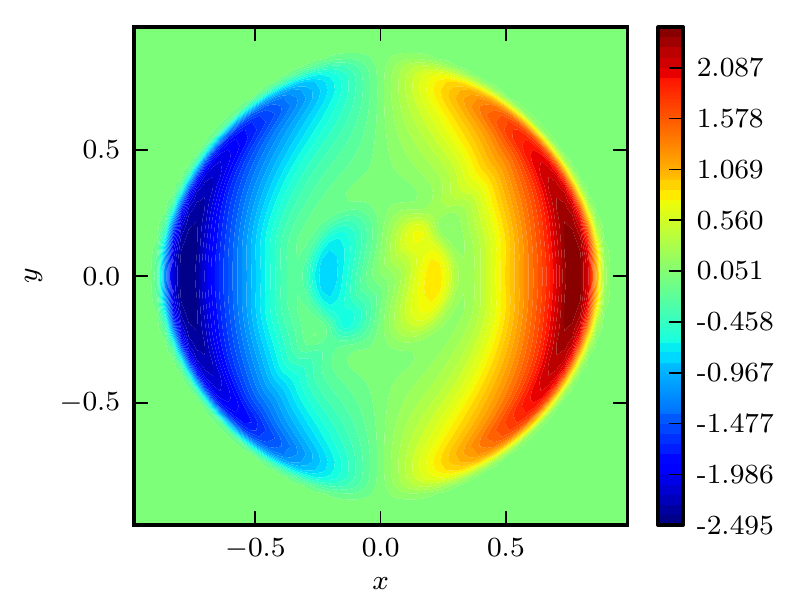} \\
		    Velocity ($y-$component)                        & Temperature                                    \\
		    \includegraphics{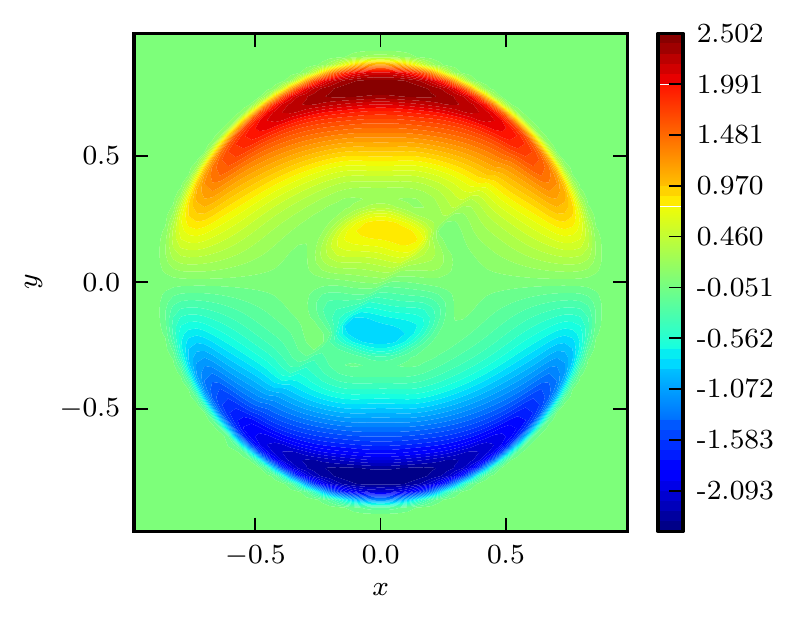}  & \includegraphics{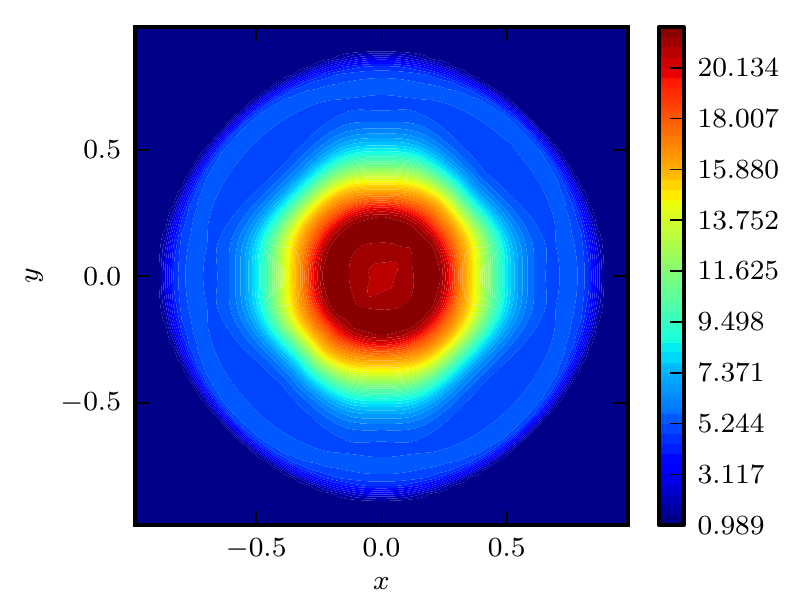}
		  \end{tabular}	  \\
		  \caption{Test 5 -- Contour plot of the first hydrodynamic quantities for the Riemann problem in temperature, given by the rescaling velocity method, at time  $t=0.1$, for $\ve = 0.05$. }
	    \label{figRiemannPbKin}
	  \end{center}
	\end{figure}

\begin{figure}
	  \begin{center}
	    Hydrodynamic model
		  \begin{tabular}{cc}
		    Density                                      & Velocity ($x-$component)                    \\
		    \includegraphics{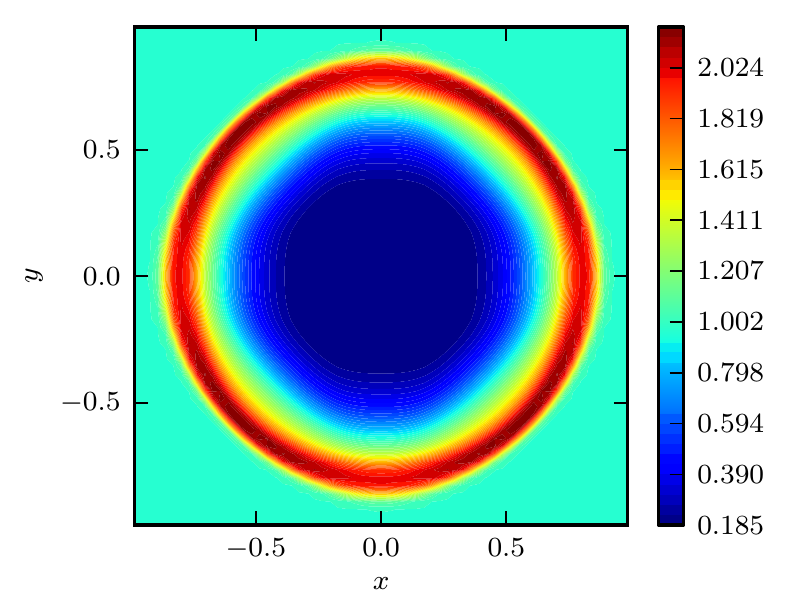} & \includegraphics{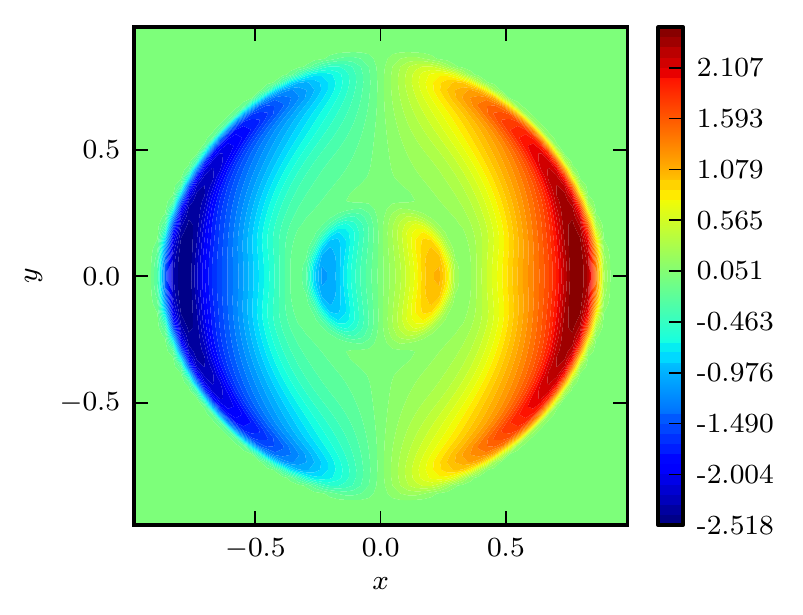} \\
		    Velocity ($y-$component)                     & Temperature                                 \\
		    \includegraphics{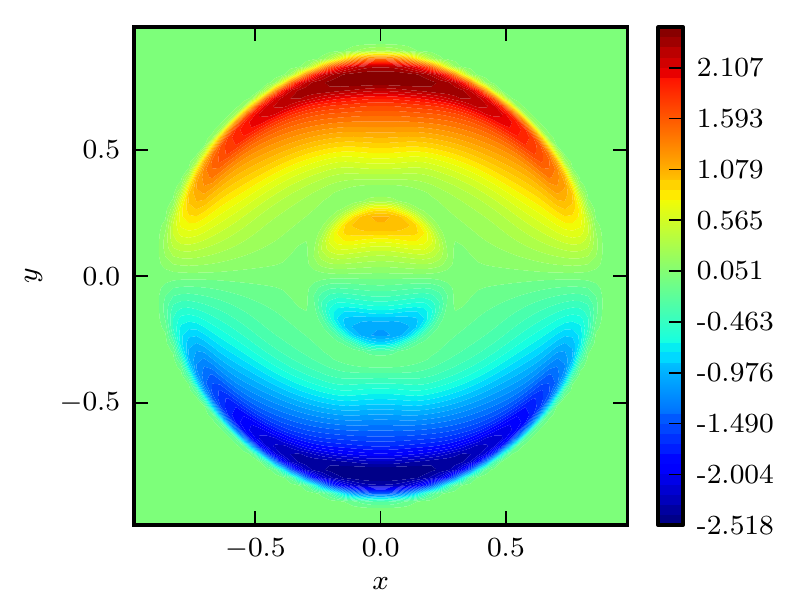}  & \includegraphics{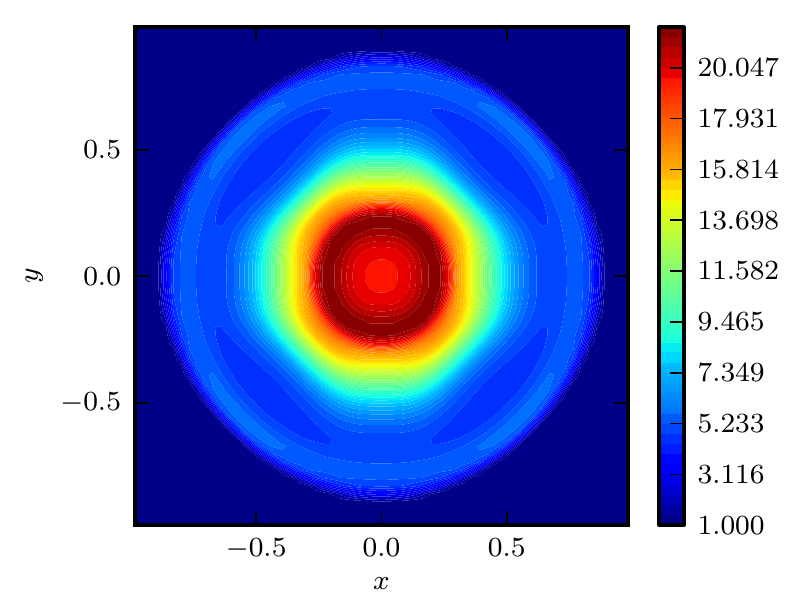}
		  \end{tabular}
		  \caption{Test 5 -- Contour plot of the first hydrodynamic quantities for the Riemann problem in temperature, given by the compressible Euler system at time  $t=0.1$. }
	    \label{figRiemannPb}
	  \end{center}
	\end{figure}

	  We present in Figures \ref{figRiemannPbKin} and \ref{figRiemannPb} the contour plot of the density, the $x$ and $y$ component of the velocity, and the temperature, at time $t=0.1$ (such a short time is needed because the front is evolving fast, due to the very large shock), for respectively the rescaling velocity method and the Euler system.
	  We observe that our method is able to compute the evolution of this problem in a sparse velocity grid, which is not the case of the solver in classical variables, and that our results are in good agreement with the hydrodynamic ones.
	  More particularly, the temperature profile are very close, which is one of the reasons that the rescaling method is working well.
    We just observe some slight oscillations in the kinetic case, but nothing really bad (and that can be overcome at a greater computational cost by using a more refined grid).

\section*{Acknowledgment}

  The authors would like to thank ICERM at Brown University for the kind hospitality during the 2011 Fall Semester Program on "Kinetic Theory and Computation" where part of this work has been carried on.  
  TR also like to thanks Lorenzo Pareschi for interesting discussions about the BGK-like simplified granular gases model.
  This work is partially supported by the European Research Council ERC Starting Grant 2009, project 239983-NuSiKiMo, by NSF grants  \#1008397 and \#1107444 (KI-Net) and ONR grant \#000141210318.
			
\bibliographystyle{acm}
\bibliography{biblioFR}

\begin{thebibliography}{10}

\bibitem{AndriesPerthame:2000}
{\sc Andries, P., {Le Tallec}, P., Perlat, J.-p., and Perthame, B.}
\newblock {The Gaussian-BGK model of Boltzmann equation with small Prandtl
  number}.
\newblock {\em European Journal of Mechanics - B/Fluids 19}, 6 (Nov. 2000),
  813--830.

\bibitem{Astillero:2004}
{\sc Astillero, A., and Santos, A.}
\newblock {A granular fluid modeled as a driven system of elastic hard
  spheres}.
\newblock In {\em The Physics of Complex Systems. New Advances and
  Perspectives}, F.~Mallamace and H.~E. Stanley, Eds., no.~2. IOS Press,
  Amsterdam, Sept. 2004, pp.~475--480.

\bibitem{Astillero:2005}
{\sc Astillero, A., and Santos, A.}
\newblock {Uniform shear flow in dissipative gases: Computer simulations of
  inelastic hard spheres and frictional elastic hard spheres}.
\newblock {\em Phys. Rev. E 72}, 3 (Sept. 2005), 1--23.

\bibitem{Besse:2011}
{\sc Besse, C., Borghol, S., Goudon, T., Lacroix-Violet, I., and Dudon, J.-P.}
\newblock {Hydrodynamic Regimes, Knudsen Layer, Numerical Schemes: Definition
  of Boundary Fluxes}.
\newblock {\em Adv. Appl. Math. Mech. 3}, 5 (2011), 519--561.

\bibitem{Bird:1994}
{\sc Bird, G.}
\newblock {\em {Molecular Gas Dynamics and the Direct Simulation of Gas
  Flows}}, 2nd~ed.
\newblock Oxford University Press, 1994.

\bibitem{Bisi:2010}
{\sc Bisi, M., Carrillo, J., and Spiga, G.}
\newblock {Some alternative methods for hydrodynamic closures to dissipative
  kinetic models}.
\newblock {\em Eur. Phys. J. Special Topics 179}, 1 (May 2010), 165--178.

\bibitem{bobylev:2000}
{\sc Bobylev, A., Carrillo, J., and Gamba, I.}
\newblock {On some properties of kinetic and hydrodynamic equations for
  inelastic interactions}.
\newblock {\em J. Stat. Phys. 98}, 3 (2000), 743--773.

\bibitem{Bobylev:1997}
{\sc Bobylev, A., and Rjasanow, S.}
\newblock {Difference scheme for the Boltzmann equation based on the fast
  Fourier transform}.
\newblock {\em Eur. J. Mech. B Fluids 16}, 2 (1997), 293--306.

\bibitem{BoRj:99}
{\sc Bobylev, A.~V., and Rjasanow, S.}
\newblock Fast deterministic method of solving the {B}oltzmann equation for
  hard spheres.
\newblock {\em Eur. J. Mech. B Fluids 18}, 5 (1999), 869--887.

\bibitem{BoRj:98}
{\sc Bobylev, A.~V., and Rjasanow, S.}
\newblock Numerical solution of the boltzmann equation using a fully
  conservative difference scheme based on the fast fourier transform.
\newblock {\em Transport Theory Statist. Phys. 29}, 3-5 (2000), 289--310.

\bibitem{brilliantov:2004}
{\sc Brilliantov, N., and P\"oschel, T.}
\newblock {\em {Kinetic {T}heory of {G}ranular {G}ases}}.
\newblock Oxford University Press, USA, 2004.

\bibitem{canuto:88}
{\sc Canuto, C., Hussaini, M., Quarteroni, A., and Zang, T.}
\newblock {\em Spectral methods in fluid dynamics}.
\newblock Springer Series in Computational Physics. Springer-Verlag, New York,
  1988.

\bibitem{Carrillo:2006Self}
{\sc Carrillo, J., Francesco, M.~D., and Toscani, G.}
\newblock {Intermediate Asymptotics Beyond Homogeneity and Self-Similarity:
  Long Time Behavior for $u_t = \Delta \phi (u)$}.
\newblock {\em Arch. Ration. Mech. Anal. 180}, 1 (Jan. 2006), 127--149.

\bibitem{carrillo:2008granular}
{\sc Carrillo, J., P\"oschel, T., and Salue{\~{n}}a, C.}
\newblock {Granular hydrodynamics and pattern formation in vertically
  oscillated granular disk layers}.
\newblock {\em J. Fluid Mech. 597\/} (2008), 119--144.

\bibitem{Carrillo:2000}
{\sc Carrillo, J., and Toscani, G.}
\newblock Asymptotic {$L^1$}-decay of solutions of the porous medium equation
  to self-similarity.
\newblock {\em Indiana Univ. Math. J. 49}, 1 (2000), 113--142.

\bibitem{Carrillo:2007a}
{\sc Carrillo, J.~A., and V\'{a}zquez, J.}
\newblock {Asymptotic Complexity in Filtration Equations}.
\newblock {\em J. Evol. Equ. 7}, 3 (Feb. 2007), 471--495.

\bibitem{Cercignani:1988}
{\sc Cercignani, C.}
\newblock {\em {The Boltzmann Equation and its Applications}}.
\newblock Springer, 1988.

\bibitem{CIP:94}
{\sc Cercignani, C., Illner, R., and Pulvirenti, M.}
\newblock {\em The mathematical theory of dilute gases}, vol.~106 of {\em
  Applied Mathematical Sciences}.
\newblock Springer-Verlag, New York, 1994.

\bibitem{Coron:89}
{\sc Coron, F.}
\newblock Derivation of slip boundary conditions for the {N}avier-{S}tokes
  system from the {B}oltzmann equation.
\newblock {\em J. Statist. Phys. 54}, 3-4 (1989), 829--857.

\bibitem{desvillettes:2005}
{\sc Desvillettes, L., and Villani, C.}
\newblock {On the Trend to Global Equilibrium for Spatially Inhomogeneous
  Kinetic Systems: the Boltzmann Equation}.
\newblock {\em Invent. Math. 159}, 2 (2005), 245--316.

\bibitem{ellis:1975}
{\sc Ellis, R., and Pinsky, M.}
\newblock {The First and Second Fluid Approximations to the Linearized
  Boltzmann Equation}.
\newblock {\em J. Math. Pures Appl. 54}, 9 (1975), 125--156.

\bibitem{filbet:2011Conv}
{\sc Filbet, F., and Mouhot, C.}
\newblock {{A}nalysis of Spectral Methods for the Homogeneous {B}oltzmann
  Equation}.
\newblock {\em Trans. Amer. Math. Soc. 363\/} (2011), 1947--1980.

\bibitem{filbet:2007}
{\sc Filbet, F., Mouhot, C., and Pareschi, L.}
\newblock {Solving the Boltzmann Equation in N log2 N}.
\newblock {\em SIAM J. Sci. Comput. 28}, 3 (2007), 1029--1053.

\bibitem{filbet:2005}
{\sc Filbet, F., Pareschi, L., and Toscani, G.}
\newblock {Accurate Numerical Methods for the Collisional Motion of (Heated)
  Granular Flows}.
\newblock {\em J. Comput. Phys. 202}, 1 (2005), 216--235.

\bibitem{Filbet:2004}
{\sc Filbet, F., and Russo, G.}
\newblock {A Rescaling Velocity Method for Kinetic Equations: the Homogeneous
  Case}.
\newblock In {\em Proceedings Modelling and Numerics of Kinetic Dissipative
  Systems\/} (Lipari, 2004), Nova-Science, p.~11.

\bibitem{Fornasier:2011}
{\sc Fornasier, M., Haskovec, J., and Toscani, G.}
\newblock Fluid dynamic description of flocking via the {P}ovzner-{B}oltzmann
  equation.
\newblock {\em Phys. D. 240}, 1 (2011), 21--31.

\bibitem{gamba:2010}
{\sc Gamba, I., and Tharkabhushanam, S.}
\newblock {Shock and Boundary Structure formation by Spectral-Lagrangian
  methods for the Inhomogeneous Boltzmann Transport Equation}.
\newblock {\em J. Comput. Math. 28}, 4 (2010), 430--460.

\bibitem{Gamba:2009}
{\sc Gamba, I.~M., and Tharkabhushanam, S.~H.}
\newblock {Spectral-Lagrangian methods for collisional models of
  non-equilibrium statistical states}.
\newblock {\em J. Comput. Phys. 228}, 6 (Apr. 2009), 2012--2036.

\bibitem{Goldhirsch:1993}
{\sc Goldhirsch, I., and Zanetti, G.}
\newblock {Clustering instability in dissipative gases}.
\newblock {\em Phys. Rev. Lett. 70}, 11 (Mar. 1993), 1619--1622.

\bibitem{Golse:2005}
{\sc Golse, F.}
\newblock {The Boltzmann equation and its hydrodynamic limits}.
\newblock In {\em Handbook of Differential Equations: Evolutionary Equations
  Vol. 2}, C.~Dafermos and E.~Feireisl, Eds. North-Holland, 2005, pp.~159--303.

\bibitem{haff:1983}
{\sc Haff, P.}
\newblock {Grain flow as a fluid-mechanical phenomenon}.
\newblock {\em J. Fluid Mech. 134\/} (1983), 401--30.

\bibitem{haile:1992}
{\sc Haile, J.}
\newblock {\em Molecular dynamics simulation: elementary methods}.
\newblock Wiley professional paperback series. Wiley, 1992.

\bibitem{Jenkins:1985}
{\sc Jenkins, J., and Richman, M.}
\newblock Grad's 13-moment system for a dense gas of inelastic spheres.
\newblock {\em Arch. Ration. Mech. Anal. 87\/} (1985), 355--377.
\newblock 10.1007/BF00250919.

\bibitem{Maxwell:1867}
{\sc Maxwell, J.~C.}
\newblock On the dynamical theory of gases.
\newblock {\em Philos. Trans. Roy. Soc. London 157\/} (1867), pp. 49--88.

\bibitem{Mischler:20062}
{\sc {M}ischler, S., and {M}ouhot, C.}
\newblock {Cooling process for inelastic Boltzmann equations for hard spheres,
  Part II: Self-similar solutions and tail behavior}.
\newblock {\em J. Stat. Phys. 124}, 2 (2006), 703--746.

\bibitem{mischler:20091}
{\sc {M}ischler, S., and {M}ouhot, C.}
\newblock {S}tability, convergence to self-similarity and elastic limit for the
  {B}oltzmann equation for inelastic hard spheres.
\newblock {\em Commun. Math. Phys. 288}, 2 (2009), 431--502.

\bibitem{Naldi:2003}
{\sc Naldi, G., Pareschi, L., and Toscani, G.}
\newblock {Spectral methods for one-dimensional kinetic models of granular
  flows and numerical quasi elastic limit}.
\newblock {\em ESAIM Math. Model. Numer. Anal. 37}, 1 (Mar. 2003), 73--90.

\bibitem{pareschi:1996}
{\sc Pareschi, L., and Perthame, B.}
\newblock {A Fourier Spectral Method for Homogeneous Boltzmann Equations}.
\newblock {\em Transport Theory Statist. Phys. 25}, 3 (1996), 369--382.

\bibitem{PaRu:SINUM:2000}
{\sc Pareschi, L., and Russo, G.}
\newblock {Numerical Solution of the Boltzmann Equation I : Spectrally Accurate
  Approximation of the Collision Operator}.
\newblock {\em SIAM J. Numer. Anal. 37}, 4 (2000), 1217--1245.

\bibitem{Poschel:2005}
{\sc P\"{o}schel, T., Brilliantov, N.~V., and Schwager, T.}
\newblock {Transient clusters in granular gases}.
\newblock {\em J. Phys.: Condens. Matter 17}, 24 (June 2005), S2705--S2713.

\bibitem{Rericha:2002}
{\sc Rericha, E.~C., Bizon, C., Shattuck, M.~D., and Swinney, H.~L.}
\newblock Shocks in supersonic sand.
\newblock {\em Phys. Rev. Lett. 88\/} (Dec 2001), 014302.

\bibitem{Rey:2012EllisPinsky}
{\sc Rey, T.}
\newblock {First and Second Order Fluid Approximations to the Linearized
  Granular Gases Equation}.
\newblock In preparation.

\bibitem{Ruiz-Montero:2010}
{\sc Ruiz-Montero, M., and Brey, J.~J.}
\newblock {On the development of inhomogeneities in freely evolving granular
  gases: The shear state and beyond}.
\newblock {\em Eur. Phys. J. Special Topics 179}, 1 (May 2010), 249--262.

\bibitem{shu:1998}
{\sc Shu, C.-W.}
\newblock High order {ENO} and {WENO} schemes for computational fluid dynamics.
\newblock In {\em High-order methods for computational physics}, vol.~9 of {\em
  Lect. Notes Comput. Sci. Eng.} Springer, Berlin, 1999, pp.~439--582.

\bibitem{sone:2007}
{\sc Sone, Y.}
\newblock {\em Molecular Gas Dynamics: Theory, Techniques, And Applications}.
\newblock Modeling And Simulation in Science, Engineering And Technology.
  Birkh{\"a}user, 2007.

\bibitem{Toscani:2004}
{\sc Toscani, G.}
\newblock {Kinetic and Hydrodynamic Models of Nearly Elastic Granular Flows}.
\newblock {\em Monatsh. Math. 142}, 1 (2004), 179--192.

\bibitem{Toscani:2007}
{\sc Toscani, G.}
\newblock Hydrodynamics from the dissipative {B}oltzmann equation.
\newblock In {\em Mathematical models of granular matter}, vol.~1937 of {\em
  Lecture Notes in Math.} Springer, Berlin, 2008, pp.~59--75.

\bibitem{VanLeer:1977a}
{\sc {Van~Leer}, B.}
\newblock {Towards the ultimate conservative difference scheme III.
  Upstream-centered finite-difference schemes for ideal compressible flow}.
\newblock {\em J. Comput. Phys. 23}, 3 (Mar. 1977), 263--275.

\bibitem{Villani:2002handbook}
{\sc Villani, C.}
\newblock {\em {A review of mathematical topics in collisional kinetic
  theory}}.
\newblock Elsevier Science, 2002.

\end{thebibliography}
\end{document}